\documentclass[12pt]{amsart}

\usepackage{amsmath,amssymb,enumerate}

\usepackage{epsfig,fancyhdr,color}

\usepackage{amssymb}
\usepackage{amsmath,amsthm}
\usepackage{bbm}
\usepackage{latexsym}
\usepackage{amscd}
\usepackage{psfrag}
\usepackage{graphicx}
\usepackage[all]{xy}

\usepackage[font=small,labelfont=bf]{caption}

\usepackage{mathtools}

\usepackage{tikz-cd}

\usepackage{stackrel}

\usepackage{mathabx,epsfig}
\def\acts{\mathrel{\reflectbox{$\righttoleftarrow$}}}

\usepackage[utf8]{inputenc}

\numberwithin{equation}{section}

\definecolor{NoteColor}{rgb}{1,0,0}

\newtheorem{theorem}{\rm\bf Theorem}[section]
\newtheorem{proposition}[theorem]{\rm\bf Proposition}

\newtheorem{corollary}[theorem]{\rm\bf Corollary}
\newtheorem*{theorem 1}{\rm\bf Proposition 1}
\newtheorem*{theorem 2}{\rm\bf Proposition 2}

\theoremstyle{definition}

\theoremstyle{remark}
\newtheorem{remark}[theorem]{\rm\bf Remark}

\newtheorem{construction}[theorem]{\rm\bf Construction}

\def\interieur#1{\mathord{\mathop{\kern 0pt #1}\limits^\circ}}

\begin{document}

\title[Universal Automorphic Functions]{Sketch of a Program for Universal Automorphic Functions to Capture  Monstrous Moonshine}
\
\author{Igor Frenkel}
\address {\hskip -2.5ex Mathematics Department\\
Yale University\\
New Haven, CT 06520\\USA}
\email{frenkel-igor{\char'100}yale.edu}
\
\author{Robert Penner}
\address {\hskip -2.5ex Institut des Hautes \'Etudes Scientifiques\\
35 route des Chartres\\
Le Bois Marie\\
91440 Bures-sur-Yvette\\
France\\
{\rm and}~Mathematics Department,
UCLA\\
Los Angeles, CA 90095\\USA}
\email{rpenner{\char'100}ihes.fr}

 \date{\today}


\begin{abstract}
We review and reformulate old and prove new results about the triad $
{\rm PPSL}_2({\mathbb Z})\subseteq{\rm PPSL}_2({\mathbb R})\acts ppsl_2({\mathbb R})
$, which provides a universal generalization of the classical automorphic triad
${\rm PSL}_2({\mathbb Z})\subseteq{\rm PSL}_2({\mathbb R})\acts psl_2({\mathbb R})$.  The leading P or $p$ in the universal setting stands for {\sl piecewise}, and 
the group ${\rm PPSL}_2({\mathbb Z})$ plays at once the role of universal modular group, universal mapping class group, Thompson group $T$ and Ptolemy group.
In particular, we construct and study new framed holographic coordinates
on the universal Teichm\"uller space and its symmetry group ${\rm PPSL}_2({\mathbb R})$, the group of piecewise ${\rm PSL}_2({\mathbb R})$ homeomorphisms of
the circle with finitely many pieces, which is dense in the group of orientation-preserving homeomorphisms of the circle.  We produce a new basis of its Lie algebra $ppsl_2({\mathbb R})$ and compute the structure constants of the Lie bracket in this basis.  We define a central extension of $ppsl_2({\mathbb R})$ and compare it with the Weil-Petersson form.  Finally, we construct a ${\rm PPSL}_2({\mathbb Z})$-invariant 1-form on the universal Teichm\"uller space formally as the Maurer-Cartan form of $ppsl_2({\mathbb R})$, which suggests the full program for developing the theory of automorphic functions for the universal triad
which is analogous, as much as possible, to the classical triad.  In
the last section we discuss the representation theory of the Lie algebra $ppsl_2({\mathbb R})$ and then pursue the universal analogy for the invariant 1-form $E_2(z)dz$, which gives rise to the spin 1 representation of $psl_2({\mathbb R})$ extended by the trivial representation.  We conjecture that the corresponding automorphic representation of $ppsl_2({\mathbb R})$ yields the bosonic CFT$_2$.  Relaxing the automorphic condition from ${\rm PSL}_2({\mathbb Z})$ to its
commutant allows the increase of the space of 1-forms six-fold additively in the classical case and twelve-fold multiplicatively in our universal case.  This leads to our ultimate conjecture that we can realize the Monster CFT$_2$ via the automorphic representation for the universal triad.  This conjecture is also bolstered by the links of both the universal Teichm\"uller and the Monster CFT$_2$ theories to the three-dimensional quantum gravity.
\end{abstract}

\maketitle

\vfill

\noindent \rule{0.165\textwidth}{0.4pt}

\noindent {\tiny Keywords: universal Teichm\"uller space, Thompson group $T$, loop algebra for sl$_2$, Weil-Petersson form, Maurer-Cartan form, Monstrous Moonshine, three dimensional quantum gravity}

\eject

\section*{Introduction}
{\let\thefootnote\relax\footnote{{ \tiny It is a pleasure to acknowledge the hospitality of the Institut des Hautes
\'Etudes Scientifiques,
where this work began, and in particular the excellent assistance of Fanny
Dufour with the figures.}}}

\noindent The idea of studying a universal Teichm\"uller space that contains the union of images of all
classical Teichm\"uller spaces goes back to Lipman Bers \cite{B}, who considered the group
of all quasisymmetric mappings of the circle ${\mathbb S}^1$ to itself modulo the M\"obius group
${\rm PSL}_2({\mathbb R})$.  

\medskip

\noindent Further development of the classical Teichm\"uller spaces \cite{Pdec} and their applications to string theory led the second-named author to a new model \cite{Puniv}  of the universal
Teichm\"uller space based on the group, denoted ${\rm Homeo}_+({\mathbb S}^1)$, of all orientation-preserving self-homeomorphisms of ${\mathbb S}^1$ with the compact-open topology modulo the M\"obius group.  The problem of providing a parametrization for the new larger model of universal Teichm\"uller space and the corresponding group of homeomorphisms of ${\mathbb S}^1$ was also resolved in \cite{Puniv} by identification of the latter group with the space ${\mathcal Tess}'$ of all ideal tesselations of the hyperbolic plane together with a choice of distinguished oriented edge, or simply doe, that is,
\begin{eqnarray}\label{0.1}
{\rm Homeo}_+({\mathbb S}^1)\approx {\mathcal Tess}'.
\end{eqnarray}
In particular, the identity element in ${\rm Homeo}_+({\mathbb S}^1)$ corresponds
to a special tesselation called the Farey tesselation.  It is obtained by applying the modular
group ${\rm PSL}_2({\mathbb Z})$ to the doe running from 0 to $\infty$.  In view of the isomorphism (\ref{0.1}), one can consider coordinates on ${\mathcal Tess}
={\mathcal Tess}'/{\rm PSL}_2({\mathbb R})$ as a measure of the distortion of a given tesselation from the specified Farey tesselation by so called shearing coordinates associated to every nonoriented edge of the Farey tesselation, or equivalently
to the elements of ${\rm PSL}_2({\mathbb Z})/({\mathbb Z}/2)$, where the ${\mathbb Z}/2$ subgroup
reverses the orientation of edges.  The shearing coordinates can be realized as logarithms of cross-ratios of certain hyperbolic lengths and are invariant under the action of the M\"obius group.  As a result
we obtain an injection
\begin{eqnarray}\label{0.2}
{\rm Homeo}_+({\mathbb S}^1)/{\rm PSL}_2({\mathbb R})\to\prod_{e\in\{{\rm edges}\}}{\mathbb R}_+^e,
\end{eqnarray}
where ${\mathbb R}_+^e$ is a copy of ${\mathbb R}_+$ associated to the edge $e$.  

\medskip

\noindent Two questions arise in relation to the parametrization (\ref{0.2}): the first is how to characterize its image,
and the second is how to circumvent the factorization by the M\"obius group and obtain directly the coordinates of ${\rm Homeo}_+({\mathbb S}^1)$ itself.  To answer these questions, we introduce new holographic coordinates and framings in Sections 2 and 3 the paper.  These new coordinates are actually based upon elaborations of the shearing coordinates defined on a decorated bundle from \cite{Pdec} over the universal Teichm\"uller space which was central also in \cite{Puniv}.

\medskip

\noindent The shearing coordinates as well as our holographic coordinates admit especially
simple transformation under the dense subgroup ${\rm PPSL}_2({\mathbb R})\subseteq{\rm Homeo}_+({\mathbb S}^1)$.  This subgroup in its turn contains a discrete subgroup
${\rm PPSL}_2({\mathbb Z})\subseteq{\rm PPSL}_2({\mathbb R})$ of piecewise
${\rm PPSL}_2({\mathbb Z})$ homeomorphisms with rational breakpoints between
pieces.  Elements of ${\rm PPSL}_2({\mathbb Z})$ turn out automatically to be once-continuously
differentiable on the circle.  This pair of groups contains the classical pair
${\rm PSL}_2({\mathbb Z})\subseteq{\rm PSL}_2({\mathbb R})$, which is the first
hint towards the new extended theory of automorphic forms.

\medskip

\noindent 
Another discrete group that is indispensable in the Teichm\"uller theory is the mapping class group
or Teichm\"uller modular group.  In our universal context it is realized by the group of the Farey-type tesselations which coincide with the Farey tesselation outside of a finite polygon and known as the Ptolemy group Pt, also introduced in \cite{Puniv}.  

\medskip

\noindent More precisely,
the {\it  flip} on an edge $e$ in a tesselation is defined by
replacing it by the other diagonal of the quadrilateral complementary to 
$\cup(\tau-\{e\})$ in ${\mathbb D}$; if $e$ is the doe, then the flip on it is enhanced by
 inducing the orientation coming from the counter-clockwise rotation of $e$.
 The {\it Ptolemy group(oid)} Pt has objects given by tesselations with doe of ${\mathbb D}$ 
 which coincide with the Farey tesselation outside of a finite polygon and morphisms given by finite compositions of flips. Triangulations with doe are combinatorially rigid, and this allows flips to be labeled by edges of a fixed tesselation, so words in these labels render Pt in fact a group.  Furthermore,
 ${\rm PSL}_2({\mathbb Z})$ sits inside Pt as those tesselations which are identical to the Farey tesselation except perhaps for the location of the doe.

\medskip

\noindent The remarkable fact of the universal setting is that under the isomorphism (\ref{0.1}), these two discrete subgroups of ${\rm Homeo}_+({\mathbb S}^1)$ coincide
$
{\rm PPSL}_2({\mathbb Z})\approx {\rm Pt}.
$
Furthermore, this universal mapping class group is also isomorphic to the celebrated
Thompson group T, so
\begin{eqnarray}\label{0.3}
{\rm Thompson~T}\approx {\rm PPSL}_2({\mathbb Z})\approx {\rm Ptolemy~Pt}.
\end{eqnarray} The rich combinatorial structure of T is studied in \cite{CFP}
and numerous sequels.  In particular, T admits a presentation by
means of two generators and certain relations similar to those of the modular group.

\medskip

\noindent
The classical theory of automorphic forms on ${\rm PSL}_2({\mathbb R})$ involves,
besides the classical modular group ${\rm PSL}_2({\mathbb Z})$, also the large class of discrete subgroups $\Gamma\subseteq{\rm PSL}_2({\mathbb R})$ commensurable with the modular group.
In our universal setting, one considers a similar class of infinite discrete groups ${\rm P}\Gamma$ associated to such subgroups $\Gamma\subseteq{\rm PSL}_2({\mathbb R})$.  By definition,
${\rm P}\Gamma$ is the subgroup of ${\rm Homeo}_+({\mathbb S}^1)$ consisting of piecewise $\Gamma$ homeomorphisms with finitely many rational breakpoints.  To study these groups, one can consider again the paving, i.e., decomposition into finite-sided congruent ideal polygons, determined by the action of $\Gamma$ on the doe from 0 to $\infty$.  In Section 1, we consider the special case when $\Gamma={\rm PSL}_2({\mathbb Z})'$ is the commutant
of the modular group.  We also review the results of \cite{L}
about the class of finitely generated such groups ${\rm P}\Gamma$, which turn out to be precisely the groups of genus zero.  This is a hint towards Monstrous Moonshine, and we conclude Section 1 with the problem of characterizing those $\Gamma\subseteq {\rm PSL}_2({\mathbb Z})$
that occur in Monstrous Moonshine in these terms.

\medskip

\noindent In the classical theory of automorphic functions, besides the pair of groups
${\rm PSL}_2({\mathbb Z})\subseteq{\rm PSL}_2({\mathbb R})$, the Lie algebra
$psl_2=psl_2({\mathbb R})$ plays a pivotal role; we employ this notation $psl_2$ in the current discussion instead of the more
standard notation $sl_2=sl_2(\mathbb R)$ simply to emphasize the relationship with the
associated Lie group.  To develop the universal counterpart of the classical theory, one requires a suitable Lie algebra for the topological group ${\rm Homeo}_+({\mathbb S}^1)$.
It has been argued in \cite{MP} that this infinite-dimensional counterpart is precisely the
algebra of piecewise $sl_2$ vector fields on the circle with finitely many pieces and rational
breakpoints between them.  This Lie algebra was denoted $psl_2$, where the 
where the {\it p}
stood for {\it piecewise}. In the present paper,
this Lie algebra will  be denoted $ppsl_2=ppsl_2({\mathbb R})$ to emphasize its relationship with the group ${\rm PPSL}_2({\mathbb R})$.  In Sections 4 and 5, we continue the study of this Lie algebra $ppsl_2$. In particular following \cite{MP}, we find a basis parametrized by the edges of the Farey tesselation, or equivalently by  ${\rm PSL}_2({\mathbb Z})/({\mathbb Z}/2)$, and explicitly derive the commutation relations in this basis.  In the next Section 6, we define the central extension of the Lie algebra $ppsl_2$ viewed as a loop algebra and compare it with the central extension
given by the universal Weil-Petersson 2-form first studied in \cite{Puniv}, which naturally extends the classical Weil-Petersson K\"ahler form computed in \cite{PWP}.

\medskip

\noindent 
All three structures of our universal triad
\begin{eqnarray}\label{0.4}
{\rm PPSL}_2({\mathbb Z})\subseteq{\rm PPSL}_2({\mathbb R})\acts ppsl_2({\mathbb R})
\end{eqnarray}
are combined in our construction of an automorphic 1-form on the universal Teichm\"uller space in Section 7.  This construction is one of the main results of the paper, and it is the first step towards
our program of developing the theory of automorphic forms for the triad
(\ref{0.4}).  
Since ${\rm PPSL}_2({\mathbb R})$ allows for
irrational breakpoints, the triad (\ref{0.4}) should more properly be written
${\rm PPSL}_2({\mathbb R})\supseteq{\rm PPSL}_2({\mathbb Z})\acts ppsl_2({\mathbb R})$.  One might simply restrict to rational breakpoints for
${\rm PPSL}_2({\mathbb R})$ as well in order to ensure (\ref{0.4}), but let us not dwell on this perhaps interesting detail.

 \medskip

\noindent
Though the proposed new theory is expected to be of a higher level of complexity than the usual theory for the classical triad $${\rm PSL}_2({\mathbb Z})\subseteq{\rm PSL}_2({\mathbb R})\acts psl_2({\mathbb R}),$$
one can pursue the analogy with the classical case whenever possible.  In particular, we argue that our automorphic 1-form on the universal Teichm\"uller space has its classical counterpart in the 1-form $E_2(z)dz$ on the hyperbolic plane, the universal cover of the classical modular curve, where $E_2(z)$ is the non-holomorphic weight two covariant Eisenstein series.  

\medskip

\noindent In Appendix B, we study its lift to an automorphic function on ${\rm PSL}_2({\mathbb R})$ and show that the resulting lift generates an indecomposable representation of $psl_2({\mathbb R})$ with a one-dimesional sub-representation.  This indecomposable representation together with its conjugate, and common one-dimensional representation, can be characterized as the harmonic subspace of automorphic functions, i.e., it is annihilated by the Laplace operator.  
Though the realization of the holomorphic and antiholomorphic discrete series of weights $4,6,8,\ldots$
by the automorphic functions is well known, the case of weight 2 has not appeared in the literature to the best of our knowledge.  Yet it is exactly an analogue of this special automorphic form which arises in our construction of a ${\rm PPSL}_2({\mathbb Z})$-invariant 1-form on universal Teichm\"uller space in the formal guise of the Maurer-Cartan form of ${\rm PPSL}_2({\mathbb R})$.

\medskip

\noindent
Our example of this universal automorphic form opens a new program of research in this subject for the triad (\ref{0.4}), as discussed in Section 8.  One of the first challenges is to find an analogue of the indecomposable representation generated by the lift of the 1-form $E_2(z)dz$ in our universal setting.  Our conjecture explained in Section 8 is that this is the infinite symmetric power $S^\infty V$ of the bosonic space that appears in the classical theory, so that $S^\infty V$ is the free bosonic field and is one of the simplest CFT$_2$.

\medskip

\noindent The theory of automorphic forms instantly enlarges various spaces by relaxing the the automorphic property to smaller discrete groups.  Thus, the passage from ${\rm PSL}_2({\mathbb Z})$ to its commutant ${\rm PSL}_2({\mathbb Z})'$ of index six increases the space $V$ additively 6-fold.  Similarly, we expect that replacing ${\rm PPSL}_2({\mathbb Z})$ by
${\rm PPSL}_2({\mathbb Z})'$ will yield a correspondingly  multiplicative increase of 
$S^\infty V$.  Thus, we expect that this mechanism will allow the construction of the Monstrous
CFT$_2$ of \cite{FLM}, and in other words capture the Monster as in the title of this paper.

\medskip

\noindent
In fact at the end of Section 8, we explain this from the other end of the theory by recalling that Monstrous Moonshine allows one to attach a genus zero subgroup $\Gamma\subseteq{\rm PSL}_2({\mathbb R})$ to every conjugacy class of the Monster group.  Then one can construct as above a $\Gamma$-invariant 1-form $E^\Gamma_2(z)dz$, where $E^\Gamma_2(z)$ is again the non-holomorphic weight two covariant Eisenstein series.  This can be lifted to a representation $V^\Gamma$ of ${\rm PSL}_2({\mathbb R})$ isomorphic to $V$ but realized by a different space of automorphic functions.  Then the universal analogue is conjectured to yield a twisted Monster representation $\Lambda^\Gamma$ as in in \cite{DF}, and the comparison of $\Lambda$ and 
$\Lambda^\Gamma$ should reveal the Monster in the universal automorphic realization.  It will be interesting to see the correspondence between the Thompson-like groups ${\rm P}\Gamma\subseteq {\rm P}{\rm PSL}_2({\mathbb R})$ introduced in Section 1 and the Monster Moonshine groups $\Gamma\subseteq {\rm PSL}_2({\mathbb R})$, both of which are genus zero and satisfy certain additional properties \cite{171,CN}

\medskip

\noindent 
Our conjecture on the relationship between automorphic forms on the universal Teichm\"uller space and the Monster CFT$_2$ is strongly supported by the links of both subjects to yet another: three-dimensional quantum gravity.  In fact, the link with universal Teichm\"uller theory has a long history in the physics literature going back to \cite{V}.  More recently, a rigorous definition of the universal phase space of AdS$_3$ gravity was given in \cite{SK}, and it was
proven there that it can be identified with the cotangent bundle over the universal Teichm\"uller space in its Bers formulation.  Thus, one could expect that the space of states of AdS$_3$ quantum gravity can be realized in some class of functions on the universal Teichm\"uller space.  However, it was also understood in the physics literature that the Bers version of universal Teichm\"uller space must be enlarged to account for all states related to black holes.  Indeed, such enlargement was an original motivation for \cite{Puniv}.

\medskip

\noindent
It was argued in \cite{W} that the general philosophy of AdS$_3$/CFT$_2$ correspondence suggests that the space of states of the simplest pure quantum gravity is precisely the Monster CFT$_2$ constructed in \cite{FLM}.  This correspondence was further supported by the explanation of Monstrous Moonshine, and in particular by the mysterious genus zero property,
from the point of view of twisted state sums of three-dimensional quantum gravity \cite{DF}.
However in spite of all these tantalizing observations, a rigorous mathematical theory of three-dimensional quantum gravity is still missing. Thus, we believe that the development of the theory of universal automorphic representations can help to build a Monstertrap for all the Moonshine
and at the same time help us to understand the true nature of three-dimesional quantum gravity.




\bigskip

\section{Farey tesselation and modular group}

\bigskip

\noindent Let ${\mathbb Z}\subseteq{\mathbb Q}\subseteq{\mathbb R}\subseteq{\mathbb C}$ denote the integers, rational, real and complex numbers, respectively, with $\hat{\mathbb Z}\subseteq \hat {\mathbb Q}\subseteq\hat {\mathbb R}\subseteq\hat {\mathbb C}$ denoting their respective one-point compactifications by the point $\infty$ at infinity.
Set  $i=\sqrt{-1}$.  

\medskip

\noindent The open unit disk ${\mathbb D}$ in the complex plane ${\mathbb C}$ is identified with the
Poincar\'e disk model of the hyperbolic plane in the standard
way, where the boundary unit circle ${\mathbb S}^1$ is identified with the circle
at infinity.  Also consider the upper half-plane ${\mathcal U}=\{
z=x+iy\in {\mathbb C}:y> 0\}$.   The {\it Cayley transform} 
$$\begin{aligned}
C:({\mathcal U},\hat{\mathbb R})&\to ({\mathbb D},{\mathbb S}^1)\\
s&\mapsto {{s-i}\over{s+i}}\\
\end{aligned}$$
induces an isomorphism of pairs.

\bigskip

{\centerline{\epsfysize2.8in\epsffile{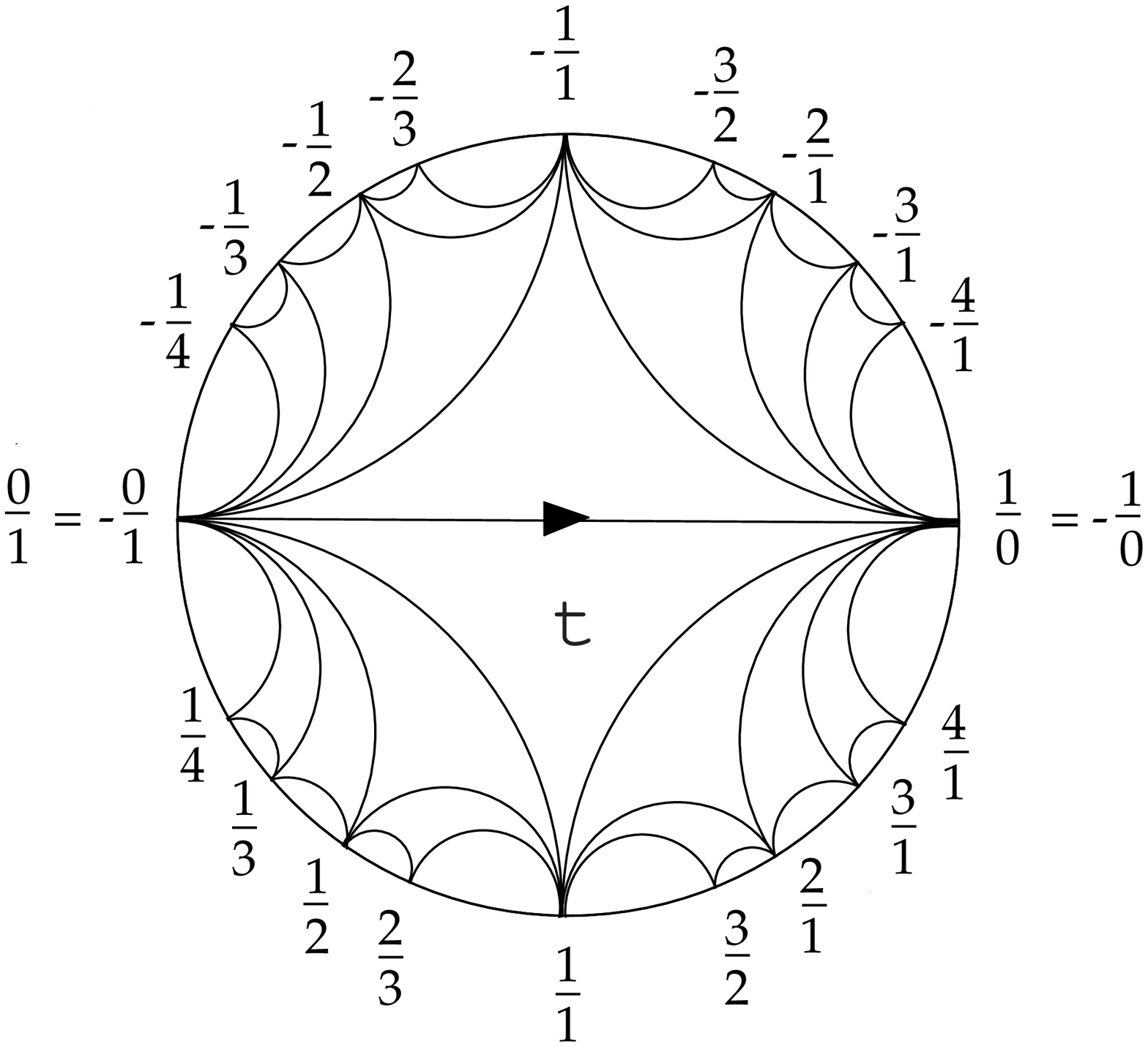}}}

{ {\bf Figure 1.} The first several generations of the Farey tesselation $\tau _*$ of the Poincar\'e disk ${\mathbb D}$
with its distinguished oriented edge.}

\bigskip

\noindent Let $t$ denote the ideal
hyperbolic triangle with vertices
$\pm1, -i\in {\mathbb S}^1$ as in Figure~1, and consider the
group ${\mathcal R}$ generated by reflections in the sides of $t$.  Define
the {\it Farey tesselation} $\tau _*$ to be the full ${\mathcal R}$-orbit
of the frontier of $t$.  $\tau_*$ has a {distinguished
oriented edge} given by the interval  from
$-1$ to $+1$.  A direct Euclidean construction of $\tau_*$
with a discussion of its history and number theoretic significance is given in \cite[$\S$3.1]{Pbook}.

\medskip

\noindent More generally, an arbitrary {\it tesselation} of ${\mathbb D}$ is a locally finite
collection $\tau$ of hyperbolic geodesics decomposing ${\mathbb D}$ into complementary
ideal triangles.  Geodesics in $\tau$
are called its {\it edges}, and $\tau$ itself is regarded
as a set of edges.  Let $\tilde\tau$ denote the set of all oriented edges of $\tau$,
and if $e\in\tilde\tau$, then let $|e|\in\tau$ denote the unoriented edge underlying $e$.  A {\it distinguished oriented edge} or {\it doe} on $\tau$ is the specification of 
an element $e\in\tilde\tau$.  Let $\tau^0\subseteq{\mathbb S}^1$ denote the set of all {\it vertices}
of $\tau$, namely, the collection of all endpoints of all edges in $\tau$.

\medskip

\noindent As is well-known, $\hat {\mathbb Q}\approx\tau _*^0$ under the 
Cayley transform, namely, 
$$C~({p/
q})={{p-iq}\over{p+iq}}={{p^2-q^2}\over{p^2+q^2}}-i~{{2pq}\over{p^2+q^2}}\in
{\mathbb S}^1,$$ as indicated in Figure~1, and we shall refer to these points as the {\it
rational points} of ${\mathbb S}^1$. We abuse notation slightly and sometimes let ${p\over
q}\in\hat{\mathbb Q}$ denote the point $C(p/q)\in {\mathbb S}^1$. 
 ${p\over q}\in\hat{\mathbb Q}$ is said to be of {\it generation g} if the radial arc in ${\mathbb D}$ from the origin to ${p\over q}$ meets the interior of $g\geq 0$ distinct ideal triangles complementary to $\tau _*$.   The standard doe of $\tau_*$ runs
 between the two rational points of generation zero, from ${0\over 1}$ to ${1\over 0}$.




\medskip

\noindent 
Another canonical tesselation of ${\mathbb D}$ with doe is the {\it dyadic tesselation} $\tau_d$, which has the same doe as the Farey tesselation, and indeed the
same generation one vertices as well, and which is recursively characterized by the property that one vertex of each triangle complementary to $\tau_d$ bisects the angle between its other two vertices.
Thus, one has $\tau_d^0=\{ e^{2\pi i(1-k)}:k\in{\mathbb Z}\}$, the points in the circle with dyadic rational arguments.  In contrast, $\tau_*^0$ consists of points in the circle with rational rectilinear coordinates.

\medskip

\noindent The {\it modular group} $${\rm PSL}_2={\rm PSL}_2({\mathbb Z})\subseteq {\mathcal R}$$ of
integral fractional linear transformations is the subgroup of ${\mathcal R}$
consisting of compositions of an even number of reflections, or in other words,
the group of two-by-two integral matrices $A$ of unit determinant modulo
the equivalence relation generated by identifying $A$ with $-A$.  More generally, the {\it M\"obius group} $${M\ddot ob}={\rm PSL}_2({\mathbb R})\supseteq {\rm PSL}_2$$
consists of the two-by-two unimodular matrices over ${\mathbb R}$ modulo the same equivalence relation.  $M\ddot ob$ is the group of orientation-preserving hyperbolic isometries of ${\mathbb D}$.

\medskip

\noindent In
particular, $A=\begin{psmallmatrix}a&b\\c&d\end{psmallmatrix}\in {\rm PSL}_2$  {\sl acts on the right}  (here following Gauss)
on the rational points by
$$A: {p\over q}~~\mapsto ~~{{pd-qb}\over{qa-pc}},$$
so the edge $e_A=(doe).A$ has initial point $-{b\over a}={{b+ia}\over{b-ia}}$ and terminal point
$-{d\over c}={{d+ic}\over{d-ic}}$.


\medskip

\noindent The {\it Thompson group} $T$ is the collection of all orientation-preserving piecewise
homeomorphisms of ${\mathbb S}^1$ with finitely many breakpoints among $\tau_d^0$ which are affine in the coordinate $\theta$ on each piece.  Recall from the Introduction the group
${\rm PPSL}_2({\mathbb Z})$ of all piecewise ${\rm PSL}_2({\mathbb Z})$ homeomorphisms of ${\mathbb S}^1$ with finitely many breakpoints among the Farey
rationals $\tau_*^0$.
In fact, the Thompson group $T$ is conjugate in Homeo$({\mathbb S}^1)_+$ to ${\rm PPSL}_2({\mathbb Z})$, where the conjugating homeomorphism fixes the endpoints of the doe and maps the Farey tesselation to the dyadic tesselation of ${\mathbb D}$ in the natural way; this conjugating homeomorphism was first studied by Minkowski in \cite{Mink} for its remarkable analytic properties.


\bigskip

\noindent Here are two standard propositions which are the starting points of our discussion:

\begin{proposition}\label{lem:farey}
The modular group ${\rm PSL}_2$ leaves setwise invariant the Farey tesselation $\tau _*$, mapping $\cup\tau _*$ onto $\cup\tau _*$.  Any 
orientation-preserving homeomorphism of the circle leaving invariant $\tau _*$ in this manner  lies in ${\rm PSL}_2$.  The modular group acts simply transitively on $\tilde\tau _*$. A fundamental domain for the action of the modular group on ${\mathcal U}$ is given by 
$\{ x+iy\in{\mathcal U}: x^2+y^2>1~{\rm and}~
|x|<{1\over 2}\}$.\hfill$\Box$
\end{proposition}

\begin{proposition}\label{lem:modular}
A generating set for ${\rm PSL}_2$ is given by any pair of
$$
~~R=\bigl (\begin{array}{cc} 0&-1\\ 1&\hskip 1.5ex 1\\\end{array}\bigr ),
~~S=\bigl (\begin{array}{cc} 0&-1\\ 1&\hskip 1.5ex 0\\\end{array}\bigr ),
~~T=\bigl (\begin{array}{cc} 1&1\\ 0&1\\\end{array}\bigr ),
~~U=\bigl (\begin{array}{cc} 1&0\\ 1&1\\\end{array}\bigr ),
$$
and $S^2=1=R^3$ is a complete set of
relations in the generators $R=T^{-1}U$ and $S=TU^{-1}T$, so ${\rm PSL}_2\approx{\mathbb Z}/2*{\mathbb Z}/3$.
In fact, $T^{-1}=R^2S$ and $U=SR^2$ are each of infinite order and are conjugate in
${\rm PSL}_2$.\hfill$\Box$
\end{proposition}

\noindent Complete relations in the generators $U$ and $T$ are given
by $ T^{-1} U T^{-1}=T U^{-1} T=U T^{-1} U=U^{-1} T U^{-1}$,
with these so-called braiding relations reflecting
the fact that ${\rm PSL}_2$ is also the mapping class group of the once-punctured
torus with $U$ and $T$ representing Dehn twists.

\medskip

\noindent Geometrically, the elliptic element
$S$ setwise fixes $|e_I|$ and reverses the orientation of the doe $e_I$,
$R$ is the elliptic transformation cyclically permuting the vertices
of the triangle to the
right of the doe,
and $U$ (respectively  $T$) is 
the parabolic transformation
with the fixed point ${0\over 1}$ (respectively ${1\over 0}$)
which cyclically permutes the incident edges of $\tau _*$ in
the counter-clockwise sense about ${0\over 1}$ (respectively
the clockwise sense about ${1\over 0}$).  Typical aspects of our enumeration
of oriented edges by elements of ${\rm PSL}_2$ are illustrated in Figure~2. 

\bigskip

{\centerline {\epsfysize2.0in\epsffile{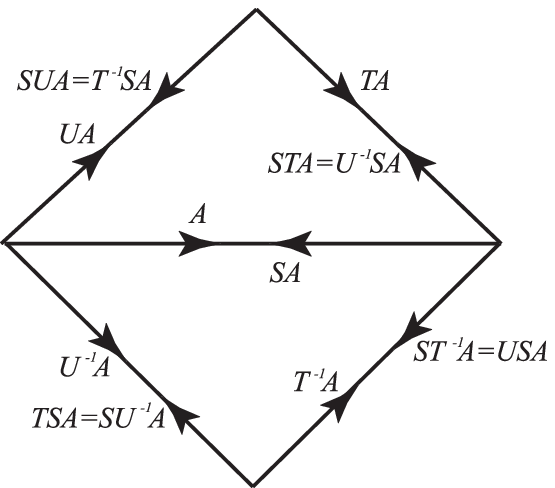}}}

\centerline{ {\bf Figure 2.} Enumeration of oriented edges near $|e_A|\in \tau_*$.}

\bigskip\bigskip

\noindent We turn our attention now to the {\it commutant}
${\rm PSL}_2'$ of ${\rm PSL}_2$, namely, the subgroup  of ${\rm PSL}_2$ generated by the 
group commutators $aba^{-1}b^{-1}$, for $a,b\in {\rm PSL}_2$.  As follows immediately
from Proposition~\ref{lem:modular}, the quotient ${\rm PSL}_2/{\rm PSL}_2'\approx {\mathbb Z}/6$.

\begin{proposition}\label{lem:six}
Consider an element of ${\rm PSL}_2$ written as a word $W=W(U,T)$ in the generators
$U$ and $T$.  Then the coset of the element $W$ in ${\rm PSL}_2/{\rm PSL}_2'$ is given
by the residue modulo six of the total exponent sum of $W$ in the letters
$U$ and $T^{-1}$.
\end{proposition}

\begin{proof}
Each of $U$ and $T^{-1}$ is infinite order and maps to the generator of
${\mathbb Z}/6\approx {\rm PSL}_2/{\rm PSL}_2'$.  Since $U$ and $T^{-1}$ are
conjugate in ${\rm PSL}_2$, the word $W$ may be written as a product of
conjugates of $U$, which thus abelianizes to the exponent sum of
$U$ in this expression of $W$ times the generator of ${\mathbb Z}/6$.
\end{proof}

\begin{corollary}\label{prop:commutant}
Fix any triangle complementary to
$\cup\tau_*$ and consider the six possible orientations on
its frontier edges.  The labels in ${\rm PSL}_2$ of these six oriented edges 
span the six commutant cosets.  Furthermore, 
${\rm PSL}_2'$ corresponds to the collection ${\mathcal C}\subseteq\tilde\tau_*$ of oriented edges determined
by the following conditions: the doe lies in ${\mathcal C}$; suppose an oriented edge
$e\in\tilde\tau_*$ lying in ${\mathcal C}$ has ideal endpoint $v\in {\mathbb S}^1$, then every third
edge of $\tau_*$ incident on $v$ also lies in ${\mathcal C}$ with alternating orientations around $v$. 
\end{corollary}

\begin{proof}
For the first part, consideration of Figure~2 shows that for the triangle to the right of the doe, we must show
that the ${\rm PSL}_2'$-cosets of $I,S,US, T^{-1} =T^{-1},TS,U^{-1}$ are distinct.  These have
respective expressions in $U,T$ given by $I$, $U T^{-1}U$, $U^2 T^{-1}U$, $ T^{-1}$, $TU T^{-1}U$, $U^{-1}$ with respective total exponent sums in $U$ and $T^{-1}$ given by $0,3,4,1,2$, $-1\equiv 5$.
The general case follows upon conjugation.

\medskip

\noindent For the second part, we likewise argue for the doe with the general result
then following by conjugation.  For the doe we must prove that $SU^{\pm 1},SU^{\pm 2}\notin{\mathcal C}$ and $SU^{\pm 3}\in{\mathcal C}$ (since obviously
$U^{\pm1},U^{\pm 2}\notin {\mathcal C}$) and likewise for $T$ instead of $U$.  To this end,
we write $SU=U T^{-1}U^2$, $SU^{-1}=U T^{-1}$, $SU^2=U T^{-1}U^3$, $SU^{-2}=U T^{-1}U^{-1}$,
$SU^3=U T^{-1}U^4$ and $SU^{-3}=U T^{-1} U^{-2}$ with respective total exponent sums in $U$ and $T^{-1}$ given by $4,2,5,-1$,$6\equiv 0,0$ as required.
\end{proof}


\medskip

\noindent The theory of the Ptolemy-Thompson group can be compared in richness only with the modular group.  It is natural to ask what is the wider class of similar groups that have similar properties as ${\rm PPSL}_2({\mathbb Z})$.  In the case of the modular group, such a class is embodied into so-called arithmetic groups, namely, subgroups of ${\rm PSL}_2({\mathbb R})$ which are commensurable with ${\rm PSL}_2({\mathbb Z})$.  They have an important common property with the modular group: their ideal compactification points in ${\mathbb D}/\Gamma$ are again $\hat{\mathbb Q}$.  Thus, the piecewise $\Gamma$ subgroups of ${\rm PPSL}_2({\mathbb R})$, denoted by P$\Gamma$, for $\Gamma$ arithmetic, might serve as a natural generalization of ${\rm PPSL}_2({\mathbb Z})$.

\medskip

\noindent The first question one might ask about these groups P$\Gamma$ is when they are finitely generated.  The answer is contained in the dissertation of Laget \cite{L}: if and only if $\Gamma$ is of genus zero.  But this is a hint towards a possible connection with Monstrous Moonshine, where the genus zero property also plays a key role.  The groups $\Gamma$ that appear in Monstrous Moonshine also have the width 1 property, namely the Farey-type tesselation or paving obtained by the action of $\Gamma$ on the edge $(0,\infty)$ in $\mathbb D$ is 1-periodic.

\medskip 

\noindent These three conditions on $\Gamma$ of arithmeticity, genus zero property and width 1 restrict the number of possible $\Gamma$ to the large but finite number 6486 according to \cite{C}.   On the other hand, there are 194 conjugacy classes in the Monster, which give rise to groups $\Gamma$
with the same three properties as above, with some identifications that reduce the number of groups $\Gamma$ involved in Monster Moonshine to 171 \cite{171}.  The problem that we want to address in this conclusion to the section is: How to characterize the Monster Moonshine $\Gamma$ among arithmetic, genus 0, width 1 groups from the properties of the piecewise $\Gamma$ groups P$\Gamma$.

\medskip

\noindent So the Monster began to show up already in Thompson-like groups.  But to capture it, we shall need other aspects of our universal automorphic triad: The universal Teichm\"uller space and the corresponding Lie algebra which we shall review and study in the next sections.

\bigskip

\section{Framed holographic coordinates}

\bigskip

\noindent We shall regard $ppsl_2$ as an appropriate limit of copies of $sl_2$, one copy at each Farey rational point in ${\mathbb S}^1$.  To take the limit, we require suitable coordinates on ideal polygons in which to 
calculate representations, and this section is dedicated to this end.  Let us begin with
the basic ``lambda length'' coordinates (sometimes called Penner coordinates) which we recall from \cite{Pdec,Pbook}.  

\medskip

\noindent A {\it decorated ideal $n$-gon} is an ideal
polygon of $n$ sides in hyperbolic space together with $n$ horocycles, one centered at each of its ideal vertices.
One such coordinate is associated to each of the $2n-3$ edges in an ideal  triangulation of the $n$-gon including its frontier edges.  As in the Introduction, there is a basic move, called a {\it flip}, on the interior edge of such a triangulation, where one removes the edge so as to produce a complimentary ideal quadrilateral and then replaces it with the other diagonal of this quadrilateral to finally produce another ideal triangulation.

\begin{theorem}\label{thm:polylam}
Fix some $n\geq 3$ and consider a {\rm decorated} ideal $n$-gon $P$ in ${\mathbb D}$.  Suppose that the frontier edges of $P$ are labeled and choose an
ideal triangulation $\Delta$ of $P$.  Then the moduli space of
such decorated polygons up to the natural action of ${M\ddot ob}$ is given by the assignment of one real lambda length $\lambda=\sqrt{{\rm exp}\, {\ell}}$ to each unoriented edge $e\in \Delta$, where $\ell$ is the signed hyperbolic distance along $e$ between the horocycles centered at its endpoints, taken with a positive sign if and only if these horocycles are disjoint.  Moreover, the lambda lengths are governed by the Ptolemy equation
$$ef=ac+bd$$
whenever $f$ arises from a flip on $e$ in the quadrilateral bounded by $a,b,c,d$ in this counter-clockwise cyclic order, here identifying an edge with its lambda length for convenience.
\hfill$\Box$
\end{theorem}


\medskip

\noindent Beyond the lambda lengths, other basic coordinates are the {\it h-lengths}  assigned to any ideal vertex $v$ of a decorated polygon $P$ by taking the hyperbolic distance along the horocycle centered at $v$ between the incident frontier edges of $P$.  
A fundamental formula illustrated in Figure 3 
relates lambda lengths and h-lengths in a triangle,
where ${{\lambda_i}\over{\lambda_j\lambda_k}}$ is the h-length opposite the edge $\lambda_i$, for $\{i,j,k\}=\{0,1,2\}$,
and we here
and hereafter again often conflate an edge with its lambda length for convenience.  It follows that
the product of h-lengths of consecutive vertices is the reciprocal square of the lambda length of
the edge they span, i.e., $\bigl({{\lambda_i}\over{\lambda_j\lambda_k}}\bigr)\bigl({{\lambda_j}\over{\lambda_i\lambda_k}}\bigr)={1\over{\lambda_k^2}}$ , so either the triple of lambda lengths
$\lambda_i$, for $i=0,1,2$ or the triple of h-lengths ${\lambda_i\over{\lambda_j\lambda_k}}$, for $\{ i,j,k\}=\{0,1,2\}$, give coordinates on the moduli space of $M\ddot ob$-orbits of decorated ideal triangles by the previous result.

\medskip

\setcounter{figure}{2}
\captionsetup[figure]{font=small,skip=0pt}
\begin{center}
\includegraphics[trim=0 0 0 0,clip,width=.45\textwidth]{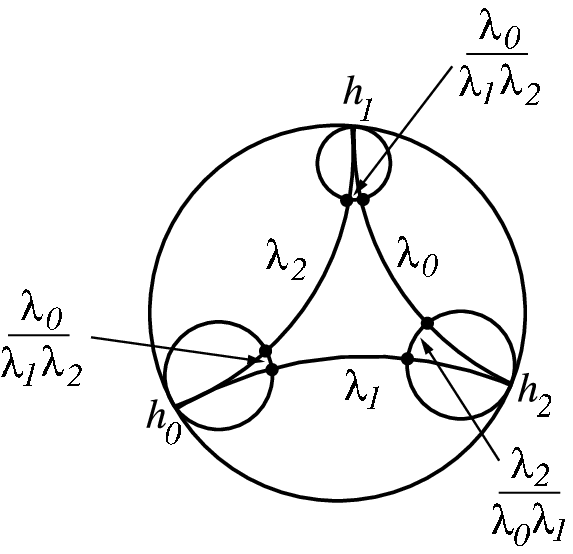}
\captionof{figure}{Three horocycles $h_i$ whose centers span a decorated ideal triangle with opposite lambda lengths $\lambda_i$
and adjacent h-lengths ${{\lambda_i}\over{\lambda_j\lambda_k}}$ for $\{i,j,k\}=\{0,1,2\}$.}
\end{center}

\bigskip

\begin{corollary} {\rm [Holographic parameters on a decorated polygon]:} \label{holoparam} For any $n\geq 3$, the moduli space
of decorated ideal $n$-gons $P$ up to the action of ${M\ddot ob}$ is parametrized
by the assignment of the lambda length of each frontier edge together with the
assignment of the h-length of each vertex of $P$.  These parameters are constrained
by three rational equations.
\end{corollary}

\begin{proof}
The proof is by induction on $n$, and the basis step $n=3$ was just discussed.  For the
induction, choose consecutive edges $a,b$ of $P$ with intermediate vertex $w$, let $T$ be the triangular convex hull of $\{a,b\}$,
$\alpha,\beta$ the respective h-lengths in $T$ opposite the endpoints of $a,b$ distinct from $w$, and let $\gamma$ denote
the h-length at $w$ of $T$.  Consider also the h-lengths $\alpha',\beta'$ of the same vertices but in $P$ rather than $T$.

\medskip

\noindent Now consider the closure $Q$ of the complement of $T$ in $P$.  Thus $Q$ is an ideal ($n-1$)-gon
which inherits all the moduli of $P$ {\sl except} that the respective h-lengths in $Q$ at $u$ and $v$ are taken to
be $\alpha'-\alpha$ and $\beta'-\beta$.  The inductive hypothesis therefore applies to $Q$ with these parameters,
and hence also to $P$ with parameters $\alpha',\beta'$ at $u,v$.  Notice the constraint that $\alpha\beta$ equals the reciprocal square of the lambda length of $T\cap Q$.  This constraint in $T$ accounts
for dimension reduction to confirm that the three rational constraints in $Q$ persist in $P$.
\end{proof}

\noindent 
The significance of this parametrization arises upon consideration of the {\it Ptolemy groupoid} {\rm Pt}$(P)$ of an ideal polygon $P$ whose objects are the ideal triangulations of $P$ with interior doe, and whose morphisms are the flips along interior edges, where the flip on the doe $e$  induces an orientation on the resulting edge and is of order four, whereas the flip on any other edge is of order two; these are the {\it face relations}.  The flips on any two edges which do not lie in the frontier of a common triangle commute, and these are the {\it commutativity relations}. The  {\it pentagon relations} arise from the serial sequence of flips alternating between two edges lying in the frontier of a common triangle, which has order five, unless one of these edges is the doe, in which case the order is ten. See \cite{Pdec,Puniv,Pbook} for further details.

\medskip

\noindent According to Proposition 4.5 in \cite{Puniv} and remarks in \cite{Pdec,Puniv}, we have

\begin{theorem}Finite sequences of flips supported on the interior of an ideal polygon $P$ act transitively on the collection of all its triangulations with interior doe.  It follows that flips generate {\rm Pt}$(P)$.  Furthermore, a complete set of relations are given by the face, commutativity  and pentagon relations.\hfill$\Box$
\end{theorem}

\medskip

\noindent So here finally is the significance of the holographic parameters: Taking the quotient by flips 
evidently renders meaningless the notion of lambda lengths of interior edges, these
edges being precisely the data are obliterated.
However, the lambda lengths of its frontier edges and the h-lengths of its vertices survive flips on interior edges to give parameters on the quotient.
This explains our mechanism of holography.

\medskip

\noindent We must go a bit further still in order to derive useful coordinates.   An ordered triple of pairwise distinct
points in the circle is a member of 
the configuration space
$$C_3=\{ (u,v,w)\in ({\mathbb S}^1)^3:u,v,w\mbox{ are pairwise distinct}\},$$
on which $M\ddot ob$ acts simply transitively
according to Proposition \ref{lem:farey}.
A {\it framing on an ideal polygon} is the specification of an arbitrary element of $C_3$, whose members need not
be among the vertices of the polygon.  A framed ideal polygon
does not necessarily contain any doe, though
an interior doe in a triangulation evidently determines an {\it associated framing} given by the
respective initial and terminal points of the doe followed by the third vertex of the triangle immediately to the right of the doe.

\medskip 

\noindent 
An ordered triple $(u,v,w)\in C_3$ uniquely determines a unit tangent vector in ${\mathbb D}$ as follows. There is a unique geodesic $g$ asymptotic to $u$ and $v$, and the orthogonal projection of $w$ onto $g$ provides a point $p\in{\mathbb D}$.  The unit tangent vector to $g$ at $p$, where $g$ is oriented with $w$ on its right, finally provides the asserted unit vector associated to $(u,v,w)$.  In these coordinates on the unit tangent bundle to ${\mathbb D}$, the almost complex structure is furthermore conveniently described by $(u,v,w)\mapsto (w,x,v)$, where ${\mathbb S}^1\ni x\neq w $ is the endpoint of the geodesic through $p$ which is asymptotic to $w$.

\begin{remark} The term ``framing'' is motivated by instantons, where the incorporation of framing
dramatically simplifies the formulas and accounts for the gauge group, cf.\ \cite{Milgram}.
The same phenomena occur here: the framing kills both the ${M\ddot ob}$ action 
and the three relations on holographic parameters
in order to produce finally useful coordinates.  
\end{remark}

\medskip

\noindent To explain the utility of the new coordinates, let us consider lambda  and h-lengths in
the upper half space model $\mathcal U$.  A decorated ideal $n$-gon $P$ is uniquely determined by a collection of pairwise disjoint points $s_i\in \hat{\mathbb R}={\mathbb R}\cup\{\infty\}$,
for $i=1,\cdots ,n$, together with a collection of Euclidean diameters $\delta_i$ of corresponding horocycles and perhaps
the $y$-coordinate $\delta_\infty$ of a horocycle about $\infty$.  Evidently these are actually coordinates, not just parameters,
as they satisfy no relations other than pairwise inequality among the $\{s_i\}_1^n$.

\medskip

\noindent 
The collection of pairs $(s_i,\delta_i)$, for $i=1,\ldots ,n$, are called {\it framed holographic coordinates}
on the space of all decorated ideal polygons, and these are our desired coordinates.  To explain the
relationship with framed polygons, take a fixed copy of ${\mathbb D}$ and choose a fixed Cayley map
$C^{-1}:\mathbb D\to\mathcal U$.  Thus, given a ${M\ddot ob}$-orbit of decorated ideal polygon in $\mathbb D$
with framing $(u,v,w)\in C_3$, there is a well-defined decorated ideal polygon in $\mathcal U$ gotten by post-composing $C$ with the unique M\"obius transformation mapping $C^{-1}(u),C^{-1}(v),C^{-1}(w)$ to the respective points ${0\over 1},{1\over 0},{1\over 1}$.

\medskip

\noindent We have proved

\begin{theorem}\label{thm:fhc} The moduli space of all framed decorated ideal $n$-gons modulo the action
of~{\rm $M\ddot ob$} is given by the pairwise distinct $n$-tuples $\{s_i\}_1^n$ of centers of horocycles in $\hat{\mathbb R}$
with Euclidean diameters $\{\delta_i\}_1^n$ called the {\rm framed holographic coordinates} $\{(s_i,\delta_i):1\leq i\leq n\}$.
The action of $A=\begin{psmallmatrix}a&b\\c&d\end{psmallmatrix}\in {M\ddot ob}$ on these coordinates is the usual right fractional linear action 
$s_i\mapsto {{ds_i-b}\over{-cs_i+a}}$ on the $\{ s_i\}_1^n$, while the $\{ \delta_i\}_1^n$ scale by the derivative 
of $A$ at $s$, so $\delta_i\mapsto{\delta_i\over{(a-cs_i)^2}}$ .
\hfill$\Box$
\end{theorem}

\begin{proof}
The only point requiring clarification is the transformation law for $\delta$ coordinates, which follows
from the fact \cite{Pdec,Pbook}
that Euclidean diameters of horocycles scale by the derivative under a M\"obius transformation.
\end{proof}

\begin{remark}
It is unsurprising that our framed holographic coordinates
are reminiscent of the dihedral coordinates of \cite{Brown}
since the latter, which pertain only to the case of planar surfaces,
are effectively related to the antecedent lambda lengths from \cite{Pdec}
for any punctured surface, which are elaborated in the new coordinates.
\end{remark}



\noindent Turning finally to stabilization given a framing
${\mathcal F}=(u,v,w)\in C_3$,
define the matrix
$$L_{\mathcal F}={1\over{(w-v)(u-w)(v-u)}}\begin{pmatrix}
v(u-w)&u(w-v)\\u-w&w-v\\
\end{pmatrix}\in {\rm SL}_2({\mathbb R}),$$
which maps $(u,v,w)\mapsto ({0\over 1},{1\over 0},{1\over 1})$ under the {right action},
and given a second framing $\bar{\mathcal F}$, define
$$\begin{aligned}
L_{\mathcal F}^{\bar{\mathcal F}}&= L_{\mathcal F}\circ L^{-1}_{\bar{\mathcal F}}=\begin{psmallmatrix}
a_{\mathcal F}^{\bar{\mathcal F}}&b_{\mathcal F}^{\bar{\mathcal F}}\\c_{\mathcal F}^{\bar{\mathcal F}}&d_{\mathcal F}^{\bar{\mathcal F}}\end{psmallmatrix}\in {\rm SL}_2({\mathbb R}).\\
\end{aligned}$$

\medskip

\noindent Now, given holographic coordinates $\{(s_i,\delta_i)\}_1^n$ with framing
${\mathcal F}=(u,v,w)$ and given $(s,\delta)$ with $s\notin\{ s_i\}_1^n$, there is a unique $i^*\in\{ 1,\ldots ,n\}$ with $s_{i^*}<s<s_{i^*+1}$.  For $j=1,\ldots ,n+1$, define
$$(\bar s_j,\bar\delta_j)=\begin{cases}
(s_j,\delta_j),&{\rm if}~ j\leq i^*;\\
(s,\delta),& {\rm if}~ j=i^*+1;\\
(s_{j+1},\delta_{j+1}),&{\rm if}~j>i^*,\\
\end{cases}$$
and given a framing $\bar{\mathcal F}=(\bar u,\bar v,\bar w)$
on $\{(\bar s_j,\bar\delta_j)\}_1^{n+1}$, finally let
$$\begin{aligned}
s_k'&= (\bar s_k) L_{\mathcal F}^{\bar{\mathcal F}},\\
\delta_k'&=(\bar\delta_k){d\over{ds}} \biggl( (\bar s_k) L_{\mathcal F}^{\bar{\mathcal F}}\biggr){},
\end{aligned}$$
for $k=1,\ldots, n+1$.

\medskip

\begin{corollary}In the notation above, the mapping
$$\{(s_i,\delta_i)\}_1^n\mapsto \{(s_i',\delta_i')\}_1^{n+1}$$
describes the stablization of framed holographic coordinates. In particular, if
the framing $(u,v,w)=({\bar u,\bar v,\bar w} )$ is constant,
then the stablization is given by inclusion and re-indexing.\hfill$\Box$
\end{corollary}

\medskip

\begin{remark}\label{rem:stablam}
Furthermore according to \cite{Pdec,Pbook}, the lambda length between 
$(s,\delta)$ and $(s',\delta')$ for $s,s'\neq\infty$ is given by
$(\delta\delta')^{-{1\over 2}}{{|s-s'|}}$
and between $(s,\delta)$ and $(\infty,\delta_\infty)$ is given by
$(\delta_\infty/\delta)^{1\over 2}$. 
On the level of the holographic parameters in Proposition \ref{holoparam}
and in the foregoing notation, let $c$ denote the lambda length of the edge
between $s_{i^*}$ and $s_{i^*+1}$,
and let $h_{i^*}$ and $h_{i^*+1}$ denote the respective nearby h-lengths.  If $a,b$ denote the respective lambda lengths of the decorated edges between
$s$ and $s_{i^*},s_{i^*+1}$, then after stabilization, the h-lengths at
$s_{i^*}$, $s$ and $s_{i^*+1}$ are respectively given by $h_{i^*}+{b\over{ac}}$, 
${c\over{ab}}$ and $h_{i^*+1}+{a\over{bc}}$.
\end{remark}

\bigskip

\section{Framed and decorated homemorphisms of ${\mathbb S}^1$}\label{sec:tess}

\bigskip

\noindent   Having developed framings on polygons as an alternative to the specification of a doe in the previous section, let us revisit from this framed point of view a basic result, Theorem 2.3 of \cite{Puniv}, whose proof we first recall.

\begin{theorem} The space
$${\mathcal Tess}'=\{{\rm tesselations~of}~{\mathbb D}~{\rm with~doe}\}$$
with topology induced by the Hausdorff topology on $\cup\tau\subset{\mathbb D}$
is naturally homeomorphic to the space $${\rm Homeo}_+\approx {\mathcal Tess}'$$
with the compact-open topology.
\end{theorem}

\begin{proof} 
Given $f\in{\rm Homeo}_+$, the image $\tau=f(\tau_*)$ of the Farey tesselation is another tesselation of ${\mathbb D}$, and the canonical doe
$e_I\in\tau_*$ maps to the doe $e=f(e_I)$ of $\tau'=(\tau,e)\in{\mathcal Tess}'$ by definition.  For the inverse map given a tesselation $\tau'=(\tau,e)$ with doe $e\in\tau$, begin by mapping the respective initial and terminal points of $e_I$ to the initial and terminal points of $e$.  Continue by mapping the further points in $\tau_*^0$ of the respective triangles containing $|e_I|$ complementary to $\tau_*$ to the left and right of $e_I$ in $\tau_*$ to those to the left and right of $e$ in $\tau$.  Proceed in this way, essentially relying on the combinatorial rigidity of a tesselation with doe, in order to define $f^0:\tau_*^0\to\tau^0$, an order-preserving injection between dense subsets by construction.  Using that a tesselation is locally-finite in ${\mathbb D}$ by definition, $f^0$ is seen to be a surjection as well, which thus interpolates a homeomorphism $f_{\tau'}\in{\rm Homeo}_+$ called the {\it characteristic map} of $\tau'$.  That the bijective assignment
$\tau'\leftrightarrow f_{\tau'}$ is bicontinuous for the stated topologies is clear, completing the proof.
\end{proof}

\noindent A {\it framing on a tesselation} $\tau$ of ${\mathbb D}$ is the specification of an element of the configuration space $C_3$ of distinct triples in ${\mathbb S}^1$, whose members are not required to lie in $\tau^0$, where the framed version of the basic result is

\begin{corollary}
The space
$${\mathcal Tess}^f=\{{\rm tesselations~of}~{\mathbb D}\}$$
is naturally homeomorphic to ${\mathcal Tess}'\approx {\rm Homeo}_+$.
\end{corollary}

\begin{proof}
As before for polygons, a doe $e$ for $\tau$ determines an associated framing,
and this gives 
an inclusion ${\mathcal Tess}'\subset{\mathcal Tess}^f$.  Conversely given a framing 
$(u,v,w)$ on the tesselation $\tau$, choose a doe $e_1$ on $\tau$ with associated framing
$(u_1, v_1,w_1)$, let $L_1=L_{u,v,w}^{u_1,v_1,w_1}\in$ $M\ddot ob$ as in the previous section, and define
$$\tau'{(\tau,u,v,w,e_1)}=L_1^{-1}\circ f_{\tau_1'} (\tau_*)\in{\mathcal Tess}'$$
with its doe induced from $e_I$ on $\tau_*$.

\medskip

\noindent We claim that $\tau_1'=\tau'{(\tau,u,v,w,e_1)}$ is independent of the choice of doe $e_1$ and to this end
choose another doe $e_2$ on $\tau$ with associated framing
$(u_2,v_2,w_2)$, let $\tau_2'=\tau'{(\tau,u,v,w,e_2)}$
and $L_2=L_{u,v,w}^{u_2,v_2,w_2}$.
The composition $(L_2^{-1}\circ f_{\tau_2'})^{-1}\circ (L_1^{-1}\circ f_{\tau_1'})$
lies in PSL$_2$ by Lemma \ref{lem:farey} since it leaves $\tau_*$ invariant, and in fact, it must be the identity since it furthermore pointwise fixes $u,v,w$. It follows that $\tau_1'$ and $\tau_2'$ have the same characteristic maps, and hence $\tau_1'=\tau_2'$ as tesselations with doe.

\medskip

\noindent The assignment $(\tau, u,v,w)\mapsto L_1^{-1}\circ f_{\tau_1'}$ is thus a well-defined inverse to the inclusion ${\mathcal Tess}'\subseteq{\mathcal Tess}^f$ induced by associated framings, which is thus a homeomorphism as required.
\end{proof}

\noindent Let us emphasize that framings thus replace distinguished oriented edges
in all regards: from determining global affine coordinates in this section, to the identification of ${\mathcal Tess}^f\approx{\mathcal Tess}'$ with ${\rm Homeo}_+$, to the isomorphisms
${\rm Pt}\approx{\rm T}\approx {\rm PPSL}_2({\mathbb Z})$ as well as their actions on tesselations.

\medskip

\noindent It is straight-forward to describe the group structure on ${\mathcal Tess}^f$ induced by composition of homeomorphisms in Homeo$_+$ but difficult to visualize on the level of tesselations, except in special cases.

\begin{proposition} Consider framings ${\mathcal F}_X$ on the Farey tesselation $\tau _*$ associated 
to does $e_X$ for $X\in$PSL$_2$.  Then $$(\tau_*,{\mathcal F}_A)\circ(\tau_*,{\mathcal F}_B)=(\tau_*,{\mathcal F}_{BA}),$$
with the analogous statement holding for any tesselation.\hfill$\Box$
\end{proposition}

\noindent Next we introduce decorated versions of the foregoing spaces and begin by defining
$$
\widetilde{{\mathcal Tess}^f}=\{\mbox{decorated and  framed tesselations of } {\mathbb D}\},
$$
where a {\it decoration} on a tesselation $\tau$ is the specification of one horocycle centered at each point of $\tau^0$, so there is a natural forgetful mapping to $\widetilde{{\mathcal Tess}^f}\to{\mathcal Tess}^f$, whose fiber can be identified with ${\mathbb R}^\omega$.

\bigskip

\noindent  Taking the quotient by the $M\ddot ob$-action on framings, we have
$$\begin{aligned}
{\mathcal Tess}^f/{M\ddot ob}&\approx{\mathcal Tess}^n=\{\mbox{tesselations }
\tau:{\scriptstyle{{ 0}\over { 1}},{1\over 1},{1\over 0}}\in\tau^0\}\, ,\\
\widetilde{{\mathcal Tess}^f}/{M\ddot ob}&\approx\widetilde{{\mathcal Tess}^n}=\{\mbox{decorated tesselations }
\tau:{\scriptstyle{{ 0}\over { 1}},{1\over 1},{1\over 0}}\in\tau^0\}\, \\
\end{aligned}$$
where the superscript $n$ stands for {\it normalized}.  Notice that the Farey tesselation $\tau_*$ is itself already normalized.

\medskip

\noindent As for the analogous elaboration of circle homemorphisms, first define spaces
$${\rm Homeo}_+^n=\{f\in{\rm Homeo}_+: f(t)=t~{\rm for} ~t={\scriptstyle{{ 0}\over { 1}},{1\over 1},{1\over 0}}\}
$$
of {\it normalized} homeomorphisms, so 
${\mathcal Tess}^n\approx{\rm Homeo}_+^n$ as ${\rm Homeo}_+^n$-spaces.  Next define the decorated version
$$\widetilde{{\rm Homeo}_+^n}=\{(\tilde f,f):f\in{\rm Homeo}_+^n~{\rm covered~by}~ \tilde f:\widetilde {{Tess}^n}\to\widetilde {{Tess}^n} \}$$
with the natural group structure, making $(\tilde f, f)\mapsto f$ a group homomorphism
$\widetilde{{\rm Homeo}_+^n}\to {\rm Homeo}_+^n$.  There is analogously an isomorphism $\widetilde{{\rm Homeo}_+^n}\approx\widetilde{{\mathcal Tess}^n}$ of
$\widetilde{{\rm Homeo}_+^n}$-spaces, where (identity,identity) corresponds to
the Farey tesselation with its canonical decoration $\tilde \tau_*$, namely, the one determined by the condition that the horocycles at endpoints of any common geodesic are taken to be osculating.

\medskip

\noindent This leads to the commutative cube

\hskip 1.5cm\begin{tikzcd}
\widetilde{{\mathcal Tess}^f} \arrow[rr,"\approx"] \arrow[dr,swap] \arrow[dd,swap] &&
  \widetilde{\rm Homeo}_+ \arrow[dd,swap] \arrow[dr] \\
& {\mathcal Tess}^f\arrow[rr,crossing over,"\approx" near start] &&
{\rm Homeo}_+\arrow[dd] \\
\widetilde{{\mathcal Tess}^n}\arrow[rr,"\approx" near end] \arrow[dr] && \widetilde{{\rm Homeo}_+^n} \arrow[dr] \\
&{\mathcal Tess}^n \arrow[rr,"\approx"] \arrow[uu,<-,crossing over]&&{\rm Homeo}_+^n
\end{tikzcd}

\medskip

\noindent where $\widetilde{{\rm Homeo}_+}$ is defined by pull-back, the vertical maps are principal
$M\ddot ob$-bundles, and the maps out of the plane of the page are ${\mathbb R}^\omega$-bundles.

\medskip

\begin{theorem}
Lambda lengths on edges of $\tau_*$ give global affine coordinates on $\widetilde{{\mathcal Tess}^n}\subseteq {\mathbb R}_+^{\tau_*}$.  Indeed, there is an explicit construction of a $M\ddot ob$-orbit of tesselations of ${\mathbb D}$ with doe from the assignment of a putative lambda length
to each edge of $\tau_*$.\hfill$\Box$
\end{theorem}

\noindent The proof of this result as Theorem 3.1 in \cite{Puniv} provides the explicit recursive construction of a $M\ddot ob$-orbit of decorated tesselations of ${\mathbb D}$ with doe from the assignment of a putative lambda length
to each edge of $\tau_*$, in analogy to the proof of Theorem \ref{thm:polylam}.  The canonical decoration $\tilde\tau_*$ on $\tau_*$, or in other words the identity element of the group $\widetilde{{\rm Homeo}_+^n}$, corresponds to taking all these lambda length coordinates equal to unity.  In the classical setting of a suitably decorated punctured surface of finite topological type uniformized by a torsion free subgroup of finite-index in $PSL_2$, or a so-called {\it punctured arithmetic surface}, the classical lambda length coordinates
on decorated Teichm\"uller space from \cite{Pdec} are likewise all unity, so this identity element in
$\widetilde{{\rm Homeo}_+^n}$ may be regarded as the universal punctured arithmetic surface.

\begin{remark}
It is not difficult to descend these lambda length coordinates on the bottom-back of the previous commutative cube to its bottom-front by assigning cross ratios ${{ac}\over{bd}}$, or shear coordinates ${\rm ln}{{ac}\over{bd}}$, instead of lambda lengths $e$ to each edge, where $e$ is the diagonal of the quadrilateral with frontier edges $a,b,c,d$ in this cyclic order; see \cite{Puniv}.
\end{remark}

\medskip

\noindent We next discuss framed tesselations as limits
of framed polygons.  To this end, notice
that there is the natural linear ordering on $\hat{\mathbb Q}$ arising from
the lexicographic ordering on pairs given by Farey generation followed by
the counter-clockwise order in ${\mathbb S}^1$ starting from ${0\over 1}$, and there is thus an induced
bijective enumeration ${\mathbb Z}_{\geq 0}\to \hat{\mathbb Q}$.

\begin{remark}
This is the raison d$'\hat{\rm e}$tre for the Farey construction.  Indeed, the mineralogist 
Farey published without proof this solution to the long-standing open problem of producing an explicit bijective enumeration of the rational numbers, and the proof was essentially immediately supplied by Cauchy.
\end{remark}

\begin{construction}\label{construction}
Suppose that $S\subseteq {\mathbb S}^1$ is a countable dense subset
enumerated by the bijection $\mu:{\mathbb Z}_{\geq 0}\to S$.  Construct a function
$f_\mu:\hat{\mathbb Q}\to{S}$ as follows.  For the basis step of our recursive construction, define 
$f_\mu({0\over 1})=\mu(0)$ and $f_\mu({1\over 0})=\mu(1)\in {\mathbb S}^1$.
Recursively suppose that
the images of the Farey points 
of generation at most $g\geq 0$ have been defined.  This collection of image points
decomposes ${\mathbb S}^1$ into $2^{g+1}$ open circular intervals.  Each such interval
with endpoints $x=f_\mu({p\over q}),y=f_\mu({r\over s})$ 
contains a point $z\in S$ of least index $\mu^{-1}(z)$,  
and we define $f_\mu({{p+r}\over{q+s}})=z$ in this case, thereby extending the function
$f_\mu$ to generation $g+1$.  We furthermore derive a collection $\tau_\mu=\{f_\mu(e):e\in\tau_*\}$
of geodesics with doe $f_\mu(e_I)$, where $f_\mu(e)$ denotes the geodesic in ${\mathbb D}$ with endpoints the $f_\mu$-image of the endpoints of $e\in\tau_*$.
\end{construction}

\noindent Suppose that $\mu:{\mathbb Z}_{\geq 0}\to S$ is a bijective enumeration
of a countable dense subset $S\subseteq{\mathbb S}^1$.  A circular interval $I\subseteq {\mathbb S}^1$ with endpoints $x,y\in S$ is said to be {\it solid} provided
$\mu^{-1}(z)> {\rm max}\{\mu^{-1}(x),\mu^{-1}(y)\}$ for every $z$ in the interior of $I$.  The enumeration $\mu$  is said to be {\it convergent} if every infinite
proper nested family $I_0\supsetneq I_1\supsetneq\cdots$ of solid open intervals
is disjoint from $S$, that is, $S\cap\bigcap_{j\geq 0} I_j=\emptyset$.

\medskip

\noindent For any $e\in\tau_*$ other than the doe, exactly one of the circular intervals in ${\mathbb S}^1$ complementary to its endpoints 
 is solid, and the Farey enumeration of $\hat{\mathbb Q}$  is convergent by definition of the Farey ordering.  We may always assume that an enumeration of a countable dense subset is indexed by the Farey tesselation in its canonical linear ordering.

\begin{proposition}
For any convergent bijective enumeration $\hat{\mathbb Q}\to S$ of a countable
dense subset $S\subset{\mathbb S}^1$, Construction \ref{construction} yields a 
bijection $f_\mu:\hat{\mathbb Q}\to S$ and a  tesselation $\tau_\mu=f_\mu(\tau_*)$ with doe.  Conversely, if $f_\mu:\hat{\mathbb Q}\to S$ is a bijection, then $\tau_\mu$ is a  tesselation, and $\mu$ must be convergent.
\end{proposition}

\begin{proof} Suppose that $\mu$ is convergent and $z\in S$.  For $g\geq 0$, let
$I_g$ denote the component of 
${\mathbb S}^1-\{ f_\mu(x):x\in\hat{\mathbb Q}\mbox{ is of generation at most } g\}$
which contains $z$.  It follows that $I_0\supsetneq I_1\supsetneq\cdots$ is a nested family of solid intervals.  This sequence must terminate, for otherwise 
$z\in S\cap\bigcap_{g\geq 0} I_g\neq\emptyset$ contradicts convergence of $\mu$.  Thus, $f_\mu$ maps onto $S$ and is injective by construction, and so is a bijection.

\medskip

\noindent To see that $\tau_\mu=f_\mu(\tau_*)$ is locally finite, suppose in order to derive a contradiction that a sequence of points in $e_i=f_\mu(e_i')$, for $i\geq 1$ and $e_i'\in\tau_*$, accumulates at some point in ${\mathbb D}$.  It follows that $e_i$ limits to some geodesic $e_\infty$ in ${\mathbb D}$.  Since the $\{e_i\}$ are pairwise disjoint in ${\mathbb D}$ by construction, we may assume that they all lie on one side of $e_\infty$.  By density of $S$ in ${\mathbb S}^1$, there is some $z\in S$ on the other side of $e_\infty$, and it lies in the $f_\mu$-image of some generation $g$ point of $\hat{\mathbb Q}$ by the already established surjectivity.  Taking a further subsequence if necessary, the complementary intervals to the endpoints of $e_i$ that contains $z$ are solid. The generation of these endpoints is therefore bounded above by $g$, so that $\{ e_i\}$ is a finite set, as required.  The same argument shows that each component of ${\mathbb D}-\cup\tau_\mu$ is an ideal triangle, so $\tau_\mu$ is indeed a tesselation.

\medskip

\noindent Conversely, suppose that $f_\mu:{\mathbb Q}\to S\subseteq {\mathbb S}^1$ is surjective and $z\in S\cap\bigcap_{j\geq 1} I_j$ for some nested sequence
$I_0\supsetneq I_1\supsetneq\cdots$ of solid intervals with endpoints in $S$.
The point $z\in S$ must a fortiori have some fixed generation, which bounds above
the generations of the endpoints of the intervals.  Since there are only finitely
many points of any given generation, the sequence must terminate, as required.
\end{proof}

\noindent It is an easy matter now to stabilize nested, decorated and framed ideal polygons:

\medskip

\begin{theorem} Consider a nested sequence $$ P_0\subseteq P_1\subseteq\cdots\subseteq  P_g\subseteq  P_{g+1}\subseteq\cdots$$
of ideal polygons with a common framing ${\mathcal F}\in C_3$, where $P_0$ is a geodesic, $P_g$ has $2^{g+1}$ sides for $g\geq 1$, and $P_{g+1}-P_g$ consists of $2^g$ ideal triangles, for all $g\geq 0$.
If the union of the ideal vertices of all the polygons is dense, then there is a well-defined
limiting tesselation also framed by ${\mathcal F}$, and any such is conversely given by such a limit.  The analogous statement holds for
decorated framed polygons and decorated framed tesselations.  Moreover, the tuple 
of framed holographic coordinates
$$(s_{0\over 1},\delta_{0\over 1}),(s_{1\over 0},\delta_{1\over 0}),
(s_{1\over 1},\delta_{1\over 1})
(s_{-{1\over 1}},\delta_{-{1\over 1}}),\ldots$$
in $(\hat{\mathbb R}\times{\mathbb R}_{>0})^\omega$ of decorated polygons in the weak topology given in their natural Farey ordering, where the $\{s_i\}$ are required to be pairwise distinct and $\cup\{s_i\}$ to be dense in $\hat{\mathbb R}$,
provide global coordinates on $\widetilde{{\mathcal Tess} ^f}\approx
\widetilde{{\rm Homeo}_+}$.\hfill$\Box$

\end{theorem}


\bigskip

\section{The Lie algebra $ppsl_2$}\label{sec:ppsl2}

\bigskip

\noindent Let $sl_2=sl_2({\mathbb R})$ denote the usual Lie algebra
of
traceless two-by-two real matrices with generators
$
e=\begin{psmallmatrix}0&1\\ 0&0 \end{psmallmatrix},
~~f=\begin{psmallmatrix}0&0\\ 1&0\end{psmallmatrix} ~{\rm and}~
~~h=\begin{psmallmatrix}+1&~~0\\ ~~~0&-1\end{psmallmatrix}$
and Lie brackets $[h,e]=2e$, $[h,f]=-2f$ and $[e,f]=h$.
(Context will distinguish the notational inconvenience of
determining between $e\in sl_2$ or $e\in\tilde\tau$.)

\medskip

\noindent Exponentiating $\begin{psmallmatrix}\alpha &\hskip1.5ex\beta \\ \gamma &-\alpha \end{psmallmatrix}
\in sl_2$
yields a one-parameter family of diffeomorphisms giving rise to a
vector field on the circle which is given by $$\{ (\gamma +\beta)~{\rm
cos}~\theta+2\alpha~{\rm sin}~\theta + (\gamma -\beta )\}
~~{\partial\over{\partial\theta}},$$
where ${\partial\over{\partial\theta}}$ denotes the constant unit vector field on ${\mathbb S}^1$
and $\theta$ its usual angular coordinate.
A vector field $A$ on ${\mathbb S}^1$ arising in this way is called a {\it
(global) $sl_2$ vector field}, and we write $A\in sl_2$ in this
case.  The values of an $sl_2$ vector field at any three distinct points in ${\mathbb S}^1$ determine it uniquely.  In case two vector fields $\vartheta_1$ and $\vartheta_2$ on ${\mathbb S}^1$ differ by a global $sl_2$ vector field, then we shall write $\vartheta_1\doteq\vartheta_2$.

\bigskip

\noindent More generally, a vector field $\vartheta$ on ${\mathbb S}^1$ is a {\it piecewise $sl_2$ vector field}
if ${\mathbb S}^1$ decomposes into {\sl finitely many} open connected circular intervals with pairwise disjoint interiors whose
endpoints are among the rational points of ${\mathbb S}^1$ so that $\vartheta$ restricts on the interior of each such interval to some global $sl_2$ vector field.  $ppsl_2$ denotes the collection of all 
such vector fields including the possibility of no breakpoints, namely $sl_2\subseteq ppsl_2$ itself.  The endpoints of the {\sl maximal such intervals} are called the {\it breakpoints} of $\vartheta$ itself.  There are no restrictions on the behavior
of $\vartheta$ at its breakpoints (except that the breakpoints, if any, are rational), and 
indeed $\vartheta\in ppsl_2$ may not even be defined at its breakpoints in general, in which case
its value is implicitly given as the average of the two one-sided limits.  

\begin{remark} Let us parenthetically recall a seminal result of Dirichlet, which he
proved in Berlin at age 24 thereby besting the top mathematicians of his day: the Fourier series of a piecewise smooth function with finitely many
pieces converges pointwise to the function itself except at the breakpoints where it is the average of its two one-sided limits.
\end{remark}

\noindent A bracket $[\vartheta_1,\vartheta_2]\in ppsl_2$ of two $\vartheta_1,\vartheta_2\in ppsl_2$
is defined in the natural way: the resulting vector field has preliminary breakpoints given by the
union of those of $\vartheta_1$ and those of $\vartheta_2$, and on each complementary
component in ${\mathbb S}^1$ of this union, the bracket is given by the usual bracket on $sl_2$.  It may happen
that the actual breakpoints of $[\vartheta_1,\vartheta_2]$ form a proper subset of
the preliminary breakpoints since we demand maximality of
intervals complementary to the actual breakpoints.

\medskip

\noindent  Consider breakpoints $\pm 1,\pm i\in {\mathbb S}^1\subseteq {\mathbb C}$, so the complementary intervals lie in 
respective quadrants I-IV in the complex plane enumerated as usual
in the counterclockwise sense beginning with quadrant I where both
coordinates are non-negative, and define
$$\Lambda(s)=\begin{cases}
\begin{psmallmatrix}s&s-s^{-1}\\0&s^{-1}\\\end{psmallmatrix},&{\rm on~quadrant~I};\\
\begin{psmallmatrix}s^{-1}&0\\s-s^{-1}&s\\\end{psmallmatrix},&{\rm on~quadrant~II};\\
\begin{psmallmatrix}s^{-1}&0\\s^{-1}-s&s\\\end{psmallmatrix},&{\rm on~quadrant~III};\\
\begin{psmallmatrix}s&s^{-1}-s\\0&s^{-1}\\\end{psmallmatrix},&{\rm on~quadrant~IV},\\
\end{cases}$$
for $s\in{\mathbb R}-\{0\}$, a one-parameter family in ${\rm PPSL}_2({\mathbb R})$.  
It is not
difficult to check that each $\Lambda(s)$ is moreover once-continuously differentiable on ${\mathbb S}^1$ including at its breakpoints.

\medskip

\noindent Since lambda lengths are $M\ddot ob$-invariant, it is clear that
$\Lambda(s)$ leaves invariant all lambda lengths on $\tau_*$ with any decoration, except the lambda length of the doe, which it scales by a factor $s$.  Thus, $\Lambda(s)$ is a one-parameter family in ${\rm PPSL}_2({\mathbb R})$ scaling just one affine lambda length coordinate, or in other words a multiplicative coordinate deformation.  

\medskip

\noindent More specifically, $\Lambda(s)\in{\rm PPSL}_2({\mathbb R})$ is a one-parameter family which on each circular interval determined by the intersection of a quadrant  with ${\mathbb S}^1$ is a hyperbolic transformation whose axis is spanned by the endpoints of the interval.  In the parlance of the Thurston school, an {\it earthquake} is an element of ${\rm PPSL}_2({\mathbb R})$ with two pieces, one of which is the identity and the other is a hyperbolic transformation whose axis is asymptotic to the breakpoints.  Thus, $\Lambda(s)$ is a special one-parameter family of compositions of four earthquakes, where the four hyperbolic translations
are chosen in order to produce homeomorphisms $\Lambda(s)$ that are once-continuously differentiable on ${\mathbb S}^1$.

\medskip

\noindent The derivative of $\Lambda(s)$ with respect to $s$ at $s=1$ is the extremely special
element $\vartheta\in ppsl_2$ illustrated in Figure~4 and called the {\it mother wavelet}.
It  is the basic building block of $ppsl_2$.  Justification for the appellate ``wavelet'' used here is given in
\cite{Philbert}.

\bigskip

{\centerline{\epsfysize2.0in\epsffile{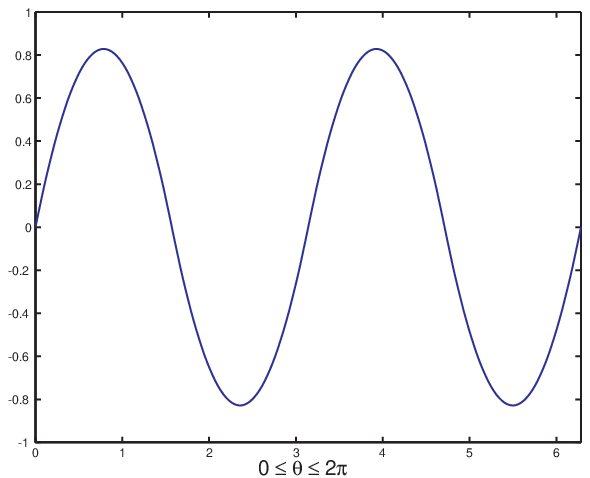}}}

\centerline{ {\bf Figure 4.} The mother wavelet $\vartheta$.}

\medskip

\noindent Direct computation confirms that the mother wavelet $\vartheta$ is given by
$$\vartheta=\begin{cases}
+h+2e,&~{\rm in~quadrant~I};\\
-h+2f,&~{\rm in~quadrant~II};\\
-h-2f,&~{\rm in~quadrant~III};\\
+h-2e,&~{\rm in~quadrant~IV}\\
\end{cases}$$

\medskip

\noindent and vanishes at each of its breakpoints, where it is once-continuously differentiable.  Notice that in general if an element of $ppsl_2$ is {\sl twice}
continuously differentiable at a breakpoint, then it is not actually a breakpoint at all, that is,
if two elements of $M\ddot ob$ agree to second order at a point, then they
must coincide.

\medskip

\noindent Let us employ the adjoint action on each piece and define the
{\it (arithmetic) wavelets}
$$\vartheta_A(\theta)=A^{-1}  \vartheta_A(\theta.A) A, ~{\rm for}~A\in{{\rm PSL}_2},$$
where the right $A$-action on $\theta\in {\mathbb S}^1$ is the natural one.
A short calculation shows that $\vartheta_S=\vartheta_I=\vartheta$, and therefore
if $A_+,A_-$ correspond to the two different orientations on a common
edge, i.e., if $A_{\mp}=SA_\pm$ then $\vartheta_{A_+}=\vartheta_{A_-}$.  
It follows that for any unoriented edge we have $|e_A|=|e_{SA}|$, and there is a corresponding
vector field $\vartheta_{|e_A|}=\vartheta_A=\vartheta_{SA}$.

\medskip

\noindent Now for each $A\in {\rm PSL}_2$, define a corresponding $X_A\in sl_2$
where $X_A$ and $\vartheta_A$ take the same values at the three points
$\pm 1,-i\in {\mathbb S}^1$, or in other words at Farey points ${0\over 1}, {1\over 0}, {1\over 1}$,
and define the {\it normalization}
$$\bar\vartheta_A=\vartheta_A-X_A, ~{\rm for}~ A\in {\rm PSL}_2.$$ 
It is not difficult to compute $X_A$ explicitly for $A\in {\rm PSL}_2$.
\bigskip

\noindent Here and herein we shall deviate slightly from the notation of previous works, which 
we first recall.

\bigskip

\noindent Following \cite{MP,Pbook} for each $A\in {\rm PSL}_2$, define for each $A\in {\rm PSL}_2$ the infinite sums
$$
\bar\phi_A=\sum _{n\geq 0} \bar\vartheta_{U^nA}~{\rm and}~~
\bar\phi_A^*=\sum _{n\leq 0} \bar\vartheta_{U^nA}
$$
respectively called {\it normalized {left} and {right fans}} and the further
infinite sums
$$
\bar\psi_A=\sum _{n\geq 0} n~\bar\vartheta_{U^nA}~{\rm and}~~
\bar\psi_A^*=\sum _{n\leq 0} n~\bar\vartheta_{U^nA}
$$
respectively called {\it normalized {left} {\rm and} {right hyperfans}}.
We have added bars to the notation for normalized (hyper)fans
from earlier work and keep the undecorated symbol for more
natural normalizations to be introduced in the next section.  It is
not difficult to compute that $\sum_{n\geq 0} X_{U^nA}$ diverges
thus explaining the need for normalization.

\bigskip

\noindent Assuming for the moment that these sums converge, as we
shall discuss presently, there is the following prescribed consequence of
this ``hyperfan formalism'' given by infinite sums of infinite sums in this way.

\medskip

\medskip

\begin{proposition}\label{lem:formal} For each $A\in {\rm PSL}_2$, we have
$$\bar\psi_A=\sum_{n\geq 1} \bar\phi_{U^nA}~{\rm and} ~~\bar\psi_A^*=\sum_{n<0} \bar\phi_{U^nA}^*,$$
as well as
$$\begin{aligned}
\bar\psi_A-\bar\psi_{UA}&=\bar\phi_{UA}~{\rm and}~~\bar\phi_A-\bar\phi_{UA}=\vartheta_A,\\
\bar\psi_{UA}^*-\bar\psi_{A}^*&=\bar\phi_{A}^*~{\rm ~~~and}~~\bar\phi_{UA}^*-\bar\phi_{A}^*=\vartheta_{A}^*\\
\end{aligned} $$
\end{proposition}
\hfill$\Box$

\medskip

\begin{corollary}\label{cor:relation}
For each $A\in {\rm PSL}_2$, equating $\bar\vartheta_A=\bar\vartheta_{SA}$ gives the so-called {\rm USA relation}
$$\begin{aligned}
\bar\psi_{UA}-2\bar\psi_A+\bar\psi_{U^{-1}A}&=\bar\psi_{USA}-2\bar\psi_{SA}+\bar\psi_{U^{-1}SA},\\
\bar\psi_{UA}^*-2\bar\psi_A^*+\bar\psi_{U^{-1}A}^*&=\bar\psi_{USA}^*-2\bar\psi_{SA}^*+\bar\psi_{U^{-1}SA}^*.\\
\end{aligned}
$$\end{corollary}
\hfill$\Box$

\medskip


\noindent The formulas in Proposition~\ref{lem:formal} and Corollary~\ref{cor:relation} thus follow immediately from the
hyperfan formalism and the symmetry $\bar\vartheta_A=\bar\vartheta_{SA}$.  The catch is showing that the
putative (hyper)fans converge.  However, notice that convergence of $\sum nx_n$ implies that of
$\sum x_n$, so the hyperfan formalism follows from convergence of hyperfans alone.

\medskip


\noindent
In the remainder of this section, we shall recall results from \cite{MP} which in particular
give sense to the {\sl normalized} (hyper)fans as elements of $ppsl_2$, and indeed  converge pointwise uniformly on compacta.
The harmonic analysis of normalized hyperfans is discussed in Appendix A.

\medskip

\noindent The {\bf first big surprise} is that in fact normalized fans and hyperfans have only 
finitely many breakpoints and lie in
$ppsl_2$, even though they are defined by infinite sums.  Indeed, a fan always
has exactly three breakpoints and is described by a continuous but not
differentiable function on the circle, and a hyperfan always has exactly
two breakpoints described by a function on the circle which is discontinuous at 
exactly one point.
Such is the nature of telescoping for normalized (hyper)fans, which are {\sl not enjoyed}
by unnormalized (hyper)fans.  

\medskip

\noindent For example,
$\bar\phi_U$ takes values $-2e$ on quadrant I, $2(h-f)$ on quadrant II
and vanishes on quadrants III and IV, while $\bar\psi_I$ takes values
$-2e$ on quadrants I and II and vanishes on quadrants III and IV.
In fact, one finds that $\bar\psi_I+\bar\psi_I^*$ is the global $sl_2$ vector field $-2e$, and
together with the fact that $\psi_A\doteq A^{-1}\psi_IA$ by construction, it follows that the span of the left hyperfans
together with $sl_2$ contains the right hyperfans.  We shall therefore henceforth restrict our attention to the former
and drop the appellate ``left'' tacitly taking only left fans and hyperfans.  There is an entirely
parallel discussion using right fans and hyperfans.

\medskip

\noindent Just as one might suspect, bracketing destroys one degree
of smoothness.  Brackets of normalized wavelets are similarly expressed as finite sums
of normalized fans while brackets of normalized fans are expressed in terms of 
normalized hyperfans.  

\medskip

\noindent The {\bf second big surprise} is that here the algebra closes with
the additive basis of normalized left hyperfans, i.e., brackets of hyperfans are
finite linear combinations of hyperfans.  

\medskip

\noindent Summarizing several of the main achievements in \cite{MP}, to which we refer the reader for proofs of the first two big surprises, we have

\begin{theorem}\label{thm:mp} 
The set of normalized left hyperfans together with the generators
$e,f,h\in sl_2$
give an additive spanning set for
the vector space $ppsl_2$ which is closed under bracketing.
Moreover, the collection of USA relations in Corollary~\ref{cor:relation},
one for each edge of $\tau_*$,
gives a complete set of relations among normalized left hyperfans. \hfill$\Box$
\end{theorem}

\noindent In fact, there is an error in the proof in \cite{MP} that normalized left
hyperfans span $ppsl_2$ which is corrected in Theorem 2.1 of \cite{Philbert} as follows.

\begin{theorem}\label{thm:basis}
Define orientations on the edges of $\tau_*$ with the orientation
from ${0\over 1}$ to ${1\over 0}$ as usual and otherwise always pointing
from lower to higher generation in the Farey enumeration, and denote
this set of oriented edges ${\mathcal O}\subset\tilde\tau_*$.  Then
${\mathcal B}=\{ \bar\psi_A:A\in {\mathcal O}\}$ together with $e,f,h\in sl_2$
is an additive basis for $ppsl_2$.\hfill$\Box$
\end{theorem}


\bigskip

\section{The new formalism for $ppsl_2$}

\bigskip

\noindent The normalization of vector fields in the previous section was necessary in order
to guarantee their pointwise convergence and assure the hyperfan formalism.  However from the point of
view of representation theory, the failure of exact equality in favor of equality $\bar\psi_A \doteq A^{-1} \bar\psi_I A$
up to global $sl_2$ vector fields introduces untoward complications.  Furthermore from the point of view of physics, the treatment of 
$sl_2$ in the previous section represents merely an additional copy of $sl_2$ essentially decoupled from the higher
Fourier modes, not the creation/annihilation and energy operators one might anticipate.

\medskip

\noindent Let us remedy both of these deficiencies by defining {\it (unnormalized) {\rm (or really differently normalized)} hyperfans}
as follows:
$$\psi_I=\begin{cases}
e,&\mbox{\rm on quadrants I and II};\\
0,&\mbox{\rm on quadrants III and IV},\\
\end{cases}$$
thus dropping the pre-factor $-2$ from before and furthermore guaranteeing the desired conjugacy formula by defining
$$\psi_A(\theta)=  A^{-1} \psi_I(\theta.A) A,  ~~{\rm for}~ A\in {\rm PSL}_2.$$
It thus follows that $B^{-1}\psi_A B=\psi_{AB}$.  

\medskip

\noindent
The {\bf third big surprise}, which is new to this paper, is that $sl_2$ is actually in the finitely supported span of these new hyperfans, as we next prove.

\begin{proposition}\label{prop:USAidentity} For $A=\begin{pmatrix}a&b\\c&d\\\end{pmatrix}\in {\rm PSL}_2$,
define 
$$\Psi_A=\psi_{STA}-2\psi_{SA}+\psi_{S T^{-1} A}-\{\psi_{UA}-2\psi_A+\psi_{ U^{-1}A}\}.$$
Then we have the identity
$$\Psi_A=\{ c(d+b)+a(d-b)\}h~~+~~(d^2-b^2+2bd)e~~+~~(a^2-c^2-2ac)f.$$
In particular
$$\Psi_I=h+e+f,~~\Psi_T=2e+f,~~\Psi_{ U^{-1}}=e+2f,$$
whence
$$\begin{pmatrix}h\cr e\cr f\end{pmatrix}={1\over 3}
\begin{pmatrix}3&-1&-1\cr 0&{~2}&-1\cr 0&-1&{~2}\cr\end{pmatrix}
\begin{pmatrix}\Psi_I\cr \Psi_T\cr \Psi_{ U^{-1}}\cr\end{pmatrix}$$
\end{proposition}

\medskip

\begin{proof} The reader will recognize $\Psi_A$ as the difference of the two sides
of the USA relation which holds for normalized hyperfans and fails by a global $sl_2$
for unnormalizaed hyperfans.  The proposition follows from direct computation
using that $\gamma-2\beta=2\alpha-\delta$, where
$$\begin{aligned}
\alpha&=  A^{-1}eA=cdh+d^2e-c^2f,~~\beta={(SA)^{-1}}eSA=abh+b^2e-a^2f,\\
\gamma&={(UA)^{-1}}e(UA)=(c+a)(d+b)h+(d+b)^2e-(c+a)^2f\\
&={(STA)^{-1}}e(STA),\\
\delta&={( U^{-1}A)^{-1}}e( U^{-1}A)=(c-a)(d-b)h+(d-b)^2e-(c-a)^2f\\
&={(S T^{-1}A)^{-1}}e(S T^{-1} A).
\end{aligned}$$
\end{proof}

\noindent We shall refer to the expression for $\Psi_A$ in terms of $e,f,h$ given in Proposition \ref{prop:USAidentity}
as the {\it USA identity} for $A\in {\rm PSL}_2$.

\begin{theorem}\label{thm:newbasis}  The hyperfans $\{ \psi_A:A\in {\rm PSL}_2\}$ span $ppsl_2$, a complete set of relations on them
is given by the USA identities as $A$ varies over ${\rm PSL}_2$, and $\{\psi_A:A\in{\mathcal O}\}$ provides a basis for $ppsl_2$, where
${\mathcal O}\subset \tilde\tau$ is given in Theorem \ref{thm:basis}.
\end{theorem}

\begin{proof}
In the old formalism, the span of normalized hyperfans together with $e,f,h$ was shown to be all of $ppsl_2$.  This contains
all unnormalized hyperfans by definition; indeed, $\psi_A\doteq\bar\psi_A$, and the reverse inclusion follows from the last part
of the previous proposition.   For the second assertion, suppose that $\sum\alpha_A\psi_A=0$ is a finite linear relation among
the unnormalized hyperfans.  It follows that then $\sum\alpha_A\bar\psi_A\doteq0$, and since both sides of this equation vanish
at ${0\over 1},{1\over 0},{1\over 1}$, it must be that in fact $\sum\alpha_A\bar\psi_A=0$, whence this relation is a consequence of the USA relations
by Theorem \ref{thm:mp}.  A non-trivial finite linear relation among $\{\psi_A:A\in{\mathcal O}\}$ likewise gives such a relation on
$\{\bar\psi_A:A\in{\mathcal O}\}$, which would contradict Theorem \ref{thm:basis}.
\end{proof}

\noindent Another favorable aspect of the new formalism is that the structure constants of $ppsl_2$ admit an explicit if not entirely trivial
expression, and it is to this end that we dedicate the remainder of this section.  The next two results follow from direct and elementary
computation, which are left to the reader.

\begin{proposition} \label{prop1} In the basis $e,f,h$ for $sl_2$, the adjoint $x\mapsto  A^{-1} x A$ of
$A=\begin{pmatrix}a&b\\ c&d\\\end{pmatrix}\in sl_2$ is given by the matrix
$$M_A=\begin{pmatrix}d^2&-b^2&2bd\\ -c^2&a^2&-2ac\\ cd&-ab&ad+bc\\
\end{pmatrix},$$
and so $$M_{A}^{-1}=M_{ A^{-1}}=\begin{pmatrix}
a^2&-b^2&-2ab\\ -c^2&d ^2&2cd\\ -ac&bd&ad+bc\\\end{pmatrix}.$$
\end{proposition}
\hfill$\Box$

\bigskip

\noindent Just as $\psi_I$ takes values $e\in\ sl_2$ on quadrants I and II and vanishes on quadrants III and IV,
we have the following analogous result for $f,h\in sl_2$.

\begin{proposition}\label{prop2}
$\psi_S+f$ and $\psi_I+\psi_{US}-\psi_S-\psi_{ U^{-1}}-f$ each have support on quadrants I and II with the former
taking there the value $f$ and the latter the value $h$.\hfill$\Box$
\end{proposition}

\noindent Now in order to compute the bracket $[\psi_B,\psi_A]$, for $A,B\in {\rm PSL}_2$, it suffices to compute
simply $[\psi_I,\psi_A]$,  since we have
$$\mbox{~~~~~~~~~~~~~~~}[\psi_B,\psi_A]~~=~~[ B^{-1}\psi_I B,\psi_A]~~=~~ B^{-1} [\psi_I,\psi_{A B^{-1}}]B.\leqno{(*)}$$
To this end with $A=\begin{pmatrix}a&b\\ c&d\\\end{pmatrix}$, there are four essential cases:

\medskip

\noindent Case 1: $0\leq -{d\over c}<-{b\over a}$, so $e_A$ lies in the bottom half plane oriented from right to left.

\medskip

\noindent Case 2: $0\leq -{b\over a}<-{d\over c}$, so $e_A$ lies in the bottom half plane oriented from left to right.

\medskip

\noindent Case 3: $ -{d\over c}<-{b\over a}\leq 0$, so $e_A$ lies in the top half plane oriented from left to right.

\medskip

\noindent Case 4: $-{b\over a}<-{d\over c}\leq 0$, so $e_A$ lies in the top half plane oriented from right to left.

\medskip

\noindent In each case, the support of $\psi_A$ lies to the left of $e_A$ by definition.  In Case 1, the supports are therefore
disjoint, and so $[\psi_I,\psi_A]=0$.  In Case 2, the support of $\psi_I$ is contained in the support of $\psi_A$, and the bracket is
supported on quadrants I and II taking value $-c^2h+2cde$, which is given according to Proposition~\ref{prop2} by
$$\begin{aligned}
-c^2(\psi_I&+\psi_{US}-\psi_S-\psi_{ U^{-1}}-f)+2cd\psi_I\\
&=c(2d-c)\psi_I+c^2(\psi_S+\psi_{ U^{-1}}-\psi_{US}+f).\\
\end{aligned}$$

\medskip

\noindent Case 3 is more challenging requiring Proposition \ref{prop1}.  We compute
$$M_{A}^{-1}\begin{pmatrix} 2cd\\ 0\\ -c^2\\
\end{pmatrix}=
\begin{pmatrix}
2ac(ad+bc)\\ -4c^3d\\ -c^2(3ad+bc)\\
\end{pmatrix},$$
whence
$$\begin{aligned}
[\psi_I,\psi_A]~=~&2ac(ad+bc)\psi_A-4c^3d(\psi_{SA}+ A^{-1} f A)\\
&-c^2(3ad+bc)\{\psi_A+\psi_{USA}-\psi_{SA}-\psi_{ U^{-1}}- A^{-1} f A\}.
\end{aligned}$$

\noindent Case 4 is still a bit more involved since $[\psi_I,\psi_A]$ perhaps has two components and can be
expressed as the difference $\mu-\nu$, where $\mu$ is supported on quadrants I and II taking value $2cde-c^2h$,
and $\nu$ is supported on the region to the left of $e_{SA}$ and taking the same value there.  Meanwhile,
Proposition \ref{prop2} gives the expression
$$\mu=c(2d-c)\psi_I+c^2(\psi_S+\psi_{ U^{-1}}-\psi_{US}+f)$$
as in Case 2.  For the other term $\nu$, we compute 
$$\begin{aligned}
M_{SA}^{-1}\begin{pmatrix}2cd\\ 0\\-c^2\\\end{pmatrix}&=\begin{pmatrix}c^2&-d^2&-2cd\\ -a^2&b^2&2ab\\ ac&-bd&-(ad+bc)\\\end{pmatrix}
\begin{pmatrix}2cd\\ 0\\-c^2\\\end{pmatrix}\\
&=\begin{pmatrix}4c^3d\\ -2ac(ad+bc)\\ c^2(3ad+bc)\\\end{pmatrix},
\end{aligned}$$
whence
$$\begin{aligned}
\nu=&4c^3d\psi_{SA}-2ac(ad+bc)(\psi_A+{(SA)}^{-1} f SA)\\
&+c^2(3ad+bc)\{ \psi_{SA}+\psi_{UA}-\psi_A-\psi_{U^{-1}SA}-{(SA)^{-1}} f SA\}.\\ 
\end{aligned}$$

\noindent Putting all this together, a little further computation finally proves

\begin{theorem}\label{thm:brackets} For $A=\begin{pmatrix} a&b\\ c&d\\\end{pmatrix}\in {\rm PSL}_2$, we have
$$[\psi_I,\psi_A]=\begin{cases}
0,~{\rm if}~0\leq -{d\over c}<-{b\over a};\\\\
c(2d-c)\psi_I+c^2(\psi_S+\psi_{ U^{-1}}-\psi_{US}+f),\\
\hskip .55cm{\rm if}~0\leq -{b\over a}<-{d\over c};\\\\
\begin{aligned}
&2ac(ad+bc)\psi_A-4c^3d(\psi_{SA}+ A^{-1} f A)\\
&-c^2(3ad+bc)\{ \psi_A+\psi_{USA}-\psi_{SA}-\psi_{ U^{-1}}- A^{-1} f A\},\\
\end{aligned}\\
\hskip .55cm{\rm if}~-{d\over c}<-{b\over a}\leq 0;\\\\
\begin{aligned}
&c(2d-c)\psi_I+c^2(\psi_S+\psi_{ U^{-1}}-\psi_{US}+f)\\
&-4c^3d\psi_{SA}+2ac(ad+bc)(\psi_A+{(SA)^{-1}} f AS)\\
&-c^2(3ad+bc)\{\psi_{SA}+\psi_{UA}-\psi_A-\psi_{ U^{-1}SA}-{(SA)^{-1}} f AS\},
\end{aligned}\\
\hskip .55cm{\rm if}~-{b\over a}<-{d\over c}\leq 0.\\
\end{cases}$$
\end{theorem}
\hfill$\Box$

\noindent It is worth noting that $[\psi_I,\psi_A]=0$ if $c=0$ and that to complete the calculation of brackets in the basis
$\{\psi_A:A\in {\mathcal O}\}$, one must still conjugate as in equation (*) and finally use the constructive proof of our Theorem \ref{thm:newbasis} given in \cite{MP} to express several of the 
resulting hyperfans in terms of this basis.  The point, however, is that brackets
are explicitly computable in the new formalism in contrast to the old.

\bigskip

\section{Three 2-forms}\label{sec:2-forms}

\bigskip

\noindent As before, let ${\partial\over{\partial\theta}}$ denote the counter-clockwise unit vector field on ${\mathbb S}^1$
in its usual angular coordinate $\theta$ and $L_n~=~ie^{in\theta} {\partial\over{\partial\theta}}$, for $n\in{\mathbb Z}$, denote the usual
generators for the {\it Witt algebra}, satisfying for $m,n\in{\mathbb Z}$ the bracket identity
$[L_m,L_n]= (m-n)~L_{m+n}$.  The Witt algebra is naturally
regarded  \cite{Segal} as the tangent space at the identity to the manifold ${\rm Diff}_+={\rm Diff}_+({\mathbb S}^1)$ of real-analytic orientation-preserving diffeomorphisms
of ${\mathbb S}^1$.  

\medskip

\noindent As the
tangent space to ${\rm Homeo}_+\supseteq{\rm Diff}_+$ as in
Section \ref{sec:tess}, $ppsl_2$ contains
the Witt algebra; this inclusion is explicitly computed in Theorem \ref{thm:Philbert}.  Conversely for $A\in{\rm PSL}_2$, the Fourier expansion of the normalized wavelets $\bar \vartheta_A$ (given in Theorem \ref{thm:fourier}) or hyperfans $\psi_A$ (given in
Proposition \ref{hyperfanfourier}), describes how $ppsl_2$ lies inside the topological closure of the Witt algebra.

\medskip

\noindent
The quotient manifold 
${\rm Diff}_+/{M\ddot ob}$, has tangent space at the identity given by
the span of $L_n$, for $n^2>1$, that is, 
the span of $e,f,h\in sl_2$ is naturally identified with the span of $L_{-1},L_0,L_{+1}$.
As likewise follows from Section \ref{sec:tess}, the quotient ${\rm Homeo}_+^n\approx{\rm Homeo_+}/{M\ddot ob}$ can be identified with the space ${\mathcal Tess}^n$ of all normalized tesselations of ${\mathbb D}$, 
and the natural bundle $\widetilde{{\mathcal Tess}^n}\to{\mathcal Tess}^n$
over this space admits global affine lambda length coordinates.

\medskip

\noindent There are the following several 2-forms defined on these spaces, which are compared in this section:

\bigskip

\noindent $\bullet$ The (pull-back from ${{\mathcal Tess}^n}$ of the) universal Weil-Petersson (WP) K\"ahler 2-form  \cite{PWP,Puniv} 
to $\widetilde{{\mathcal Tess}^n}\approx\widetilde{{\rm Homeo}_+^n}$ is given by  
$$\omega=-2\sum d{\rm log} a\wedge d{\rm log} b~~+~~
d{\rm log} b\wedge d{\rm log} c~~+~~
d{\rm log} c\wedge d{\rm log} a ,$$
where the sum is over the set of all triangles
complementary to $\tau_*$ in $\mathbb D$
and the triangle has edges, here conflated with lambda lengths as usual, in the clockwise order $a,b,c$ in the orientation on
the underlying surface.  This sum converges provided
the homeomorphism is $C^{{3\over 2}+\epsilon}$ smooth on ${\mathbb S}^1$, cf.\ \cite{Verjovsky-Nag}.

\bigskip

\noindent $\bullet$ The Kirillov-Kostant (KK) form \cite{KK, Verjovsky-Nag} is defined 
on Diff$_+$/$M\ddot ob$ by
$$\kappa_a(L_m,L_n)=a(m^3-m)\delta_{m,-n},$$
where $a\in {\mathbb C}$ and $\delta$ is the Kronecker delta function.

\bigskip 

\noindent $\bullet$ The natural loop-algebra (LA)  2-cocycle \cite{Fuchs} is defined by
$$\gamma(\psi^1,\psi^2)={1\over{2\pi}}\int_0^{2\pi} {\rm tr}\biggl ( \psi^1(\theta)\cdot d\psi^2(\theta)\biggr),$$
where $\psi_1,\psi_2\in sl_2^{{\mathbb S}^1}$ and tr denotes the trace.  Our algebra $ppsl_2$ is a sub-algebra of a certain
completion of $sl_2^{{\mathbb S}^1}$ on which this 2-cocycle still makes sense, and indeed integration by parts provides
in this case the
explicit expression 
$$
\gamma(\psi^1,\psi^2)={1\over 2} \sum_{\theta\in\Pi(\psi^2)}
{\rm tr}\biggl \{ [\psi^1(\theta^+)+\psi^1(\theta^-)]~[\psi^2(\theta^+)-\psi^2(\theta^-)]\biggr \}
$$
for the 2-cocycle $\gamma$ on $ppsl_2$,
where $\Pi(\psi)\subset\hat{\mathbb Q}\subset {\mathbb S}^1$ denotes the set of breakpoints of $\psi$
and $\theta ^\pm$ denotes a point slightly ${{\rm beyond}\over{\rm before}}$ the point $\theta\in\Pi(\psi)$
in the counter-clockwise orientation on ${\mathbb S}^1$.  One can verify directly that $\gamma$ is a 2-cocycle by
checking skew-symmetry and the 2-cocycle property, which follows from the 2-cocycle property
for the first expression in this paragraph for $\gamma$ on $sl_2^{{\mathbb S}^1}$.

\bigskip

\noindent In fact, these three 2-forms are pairwise identical up to overall constants, as we shall discuss.  The constant for WP and LA will be computed here, and the constant for KK
and WP was already calculated in Theorem 5.5 of \cite{Puniv}, namely, we have

\begin{theorem}\label{thm:KK}
 The WP K\"ahler 2-form $\omega$ and the KK form $\kappa_a$ are related by
$$\omega=\kappa_a,~~{\rm for}~a=2\pi i.$$
\end{theorem}

\noindent The computational
proof in \cite{Puniv} using the basis of normalized arithmetic wavelets is involved and delicate depending upon the constraints 
on small Fourier modes described in Remark
\ref{rmkonKK}.  

\begin{theorem}  The WP K\"ahler 2-form $\omega$ and the LA cocycle $\gamma$ are related by
$$\gamma=-4\omega.$$
\end{theorem}

\begin{proof}
Recall that $\Pi(\psi)\subset \hat {\mathbb Q}\subset {\mathbb S}^1$ denotes the (finite) set of breakpoints of $\psi\in ppsl_2$ and let
$\bar\vartheta_A=\vartheta_A-X_A$, i.e., $\bar\vartheta_A$ is the normalized arithmetic wavelet and $\vartheta_A$ the unnormalized
one.  In particular, we have
$$\Pi(\bar\vartheta_A)=\Pi(\vartheta_A)=\Pi(\vartheta).A,~~{\rm for}~A\in {\rm PSL}_2,$$
and moreover
$$\gamma(\bar\vartheta_A,\bar\vartheta_B)=\gamma(\bar\vartheta_A,\vartheta_B),~~{\rm for}~A,B\in sl_2,$$
since $X_B$ is both added and subtracted in a difference of $\bar\vartheta_B$-values
in their contribution to $\gamma$.  Thus
$$\gamma(\bar\vartheta_A,\bar\vartheta_B)=\gamma(\vartheta_A,\vartheta_B),~~{\rm for}~A,B\in sl_2$$
as well by skew symmetry of $\gamma$.

\medskip

\noindent We first claim that $\gamma(\bar\vartheta_A,\bar\vartheta_B)=0$ if $\bar\vartheta_A$ takes a common value
$Y_A\in sl_2$ at each point of $\Pi(\bar\vartheta_B)$. To see this, we compute
$$\begin{aligned}
\gamma(\bar\vartheta_A,\bar\vartheta_B)&=\gamma(\bar\vartheta_A,\vartheta_B)=\gamma(Y_A,\vartheta_B)\\
&=\sum_{\theta\in\Pi(\vartheta_B)}~~{\rm tr}\biggl \{ Y_A\cdot[\vartheta_B(\theta^+)-\vartheta_B(\theta^-)]\biggr \}  \\
&=\sum_{\theta\in\Pi(\vartheta)}~~{\rm tr}\biggl \{ Y_A\cdot[\vartheta_B(\theta^+.B)-\vartheta_B(\theta^-.B)]\biggr \}  \\
&={\rm tr} \sum_{\theta\in\Pi(\vartheta)}~~Y_A\biggl \{ B^{-1}[\vartheta(\theta^+)-\vartheta(\theta^-)]  B\biggr \}\\
&={\rm tr} \biggl \{ Y_AB^{-1}\biggl ( \sum_{\theta\in\Pi(\vartheta)}[\vartheta(\theta^+)-\vartheta(\theta^-)] \biggr ) B\biggr \}\\
&={\rm tr}\biggl\{ Y_AB^{-1}[4e+2(f-e-h)-4f+2(h+f-e)]B\biggr \}\\
&={\rm tr}~Y_AB^{-1}0B\\
&=0.\\
\end{aligned}$$

\noindent It remains to consider $\gamma(\bar\vartheta_A,\bar\vartheta_{UA})$ and $\gamma(\bar\vartheta_A,\bar\vartheta_{TA})$,
and we begin with the former.  Let $\xi_B,\eta_B$ denote the respective initial and terminal point of $e_B=e_I.B,$
where $e _I$ is the doe of $\tau_*$ as usual,
so
$$\begin{aligned}
\Pi(\bar\vartheta_A)&=\{\xi_A,\eta_A,\eta_{{U^{-1}}A},\eta_{UA}\},\\
\Pi(\bar\vartheta_{UA})&=\{\xi_A,\eta_A,\eta_{UA},\eta_{U^2A}\}.\\
\end{aligned}$$
Thus
$$\begin{aligned}
\gamma&=\gamma(\bar\vartheta_A,\bar\vartheta_{UA})=\gamma(\vartheta_A,\vartheta_{UA})\\&={1\over 2}\sum_{\theta\in\Pi(\vartheta_{UA})}
{\rm tr}\{[\vartheta_A(\theta^+)+\vartheta_A(\theta^-)][\vartheta_A(\theta^+)-\vartheta_A(\theta^-)]\},
\end{aligned}$$
so
$$\begin{aligned}
2\gamma&={\rm tr}\left\{
{\scriptstyle
\begin{aligned}
&[\vartheta_A(\xi^+.A)+\vartheta_A(\xi^-.A)]~~[\vartheta_{UA}(\xi^+.A)-\vartheta_{UA}(\xi^-.A)]\\
+&[\vartheta_A(\eta^+.A)+\vartheta_A(\eta^-.A)]~~[\vartheta_{UA}(\eta^+.A)-\vartheta_{UA}(\eta^-.A)]\\
+&[\vartheta_A(\eta^+.UA)+\vartheta_A(\eta^-.UA)]~~[\vartheta_{UA}(\eta^+.UA)-\vartheta_{UA}(\eta^-.UA)]\\
+&2[\vartheta_A(\eta^+.U^2A)]~~[\vartheta_{UA}(\eta^+.U^2A)-\vartheta_{UA}(\eta^-.U^2A)]\\
\end{aligned}
}
\right\}\\\\
&={\rm tr}\left\{
\begin{aligned}
&[\vartheta(\xi^+)+\vartheta(\xi^-)]~~[(U)^{-1}\vartheta(\xi^+.U^{-1})U-(U)^{-1}\vartheta(\xi^-.U^{-1})U]\\
+&[\vartheta(\eta^+)+\vartheta(\eta^-)]~~[(U)^{-1}\vartheta(\eta^+.U^{-1})U-(U)^{-1}\vartheta(\eta^-.U^{-1})U]\\
+&[\vartheta(\eta^+.U)+A^{-1}\vartheta(\eta^-.U)]~~[(U)^{-1}\vartheta(\eta^+)U-(U)^{-1}\vartheta(\eta^-)U]\\
+&2[\vartheta(\eta^+.U^2)]~~[(U)^{-1}\vartheta(\eta^+.U)U-(U)^{-1}\vartheta(\eta^-.U)U]\\
\end{aligned}
\right\}\\\\
&={\rm tr}\left\{
\begin{aligned}
&\biggl[\begin{psmallmatrix}-1&0\\ -2&1\\\end{psmallmatrix}+\begin{psmallmatrix}-1&0\\ \hskip1.3ex2&1\\\end{psmallmatrix}\biggr]\begin{psmallmatrix}\hskip 1.3ex1&0\\ -1&1\\\end{psmallmatrix}
\biggl [\begin{psmallmatrix}-1&0\\ -2&1\\\end{psmallmatrix}-\begin{psmallmatrix}-1&0\\\hskip1.3ex 2&1\\\end{psmallmatrix}\biggr ]\begin{psmallmatrix}1&0\\ 1&1\\\end{psmallmatrix}\\
+&\biggl[\begin{psmallmatrix}1&\hskip1.3ex2\\ 0&-1\\\end{psmallmatrix}+\begin{psmallmatrix}1&-2\\ 0&-1\\\end{psmallmatrix}\biggr]\begin{psmallmatrix}\hskip1.3ex1&0\\ -1&1\\\end{psmallmatrix}
\biggl[\begin{psmallmatrix}1&-2\\ 0&-1\\\end{psmallmatrix}-\begin{psmallmatrix}-1&0\\ -2&1\\\end{psmallmatrix}\biggr]\begin{psmallmatrix}1&0\\ 1&1\\\end{psmallmatrix}\\
+&\biggl[\begin{psmallmatrix}-1&0\\ \hskip1.3ex2&1\\\end{psmallmatrix}+\begin{psmallmatrix}1&\hskip1.3ex2\\ 0&-1\\\end{psmallmatrix}\biggr]\begin{psmallmatrix}\hskip1.3ex1&0\\ -1&1\\\end{psmallmatrix}
\biggl[\begin{psmallmatrix}1&\hskip1.3ex2\\ 0&-1\\\end{psmallmatrix}-\begin{psmallmatrix}1&-2\\ 0&-1\\\end{psmallmatrix}\biggr]\begin{psmallmatrix}1&0\\ 1&1\\\end{psmallmatrix}\\
+&2\begin{psmallmatrix}-1&0\\ \hskip1.3ex2&1\\\end{psmallmatrix}\begin{psmallmatrix}\hskip1.3ex1&0\\ -1&1\\\end{psmallmatrix}
\biggl[\begin{psmallmatrix}-1&0\\ \hskip1.3ex2&1\\\end{psmallmatrix}-\begin{psmallmatrix}1&\hskip 1.3ex2\\ 0&-1\\\end{psmallmatrix}\biggr]\begin{psmallmatrix}1&0\\ 1&1\\\end{psmallmatrix}\\
\end{aligned}
\right\}\\\\
&={\rm tr}\left\{
\begin{psmallmatrix}\hskip1.3ex0&0\\ -8&0\\\end{psmallmatrix}+\begin{psmallmatrix}\hskip1.3ex0&-4\\ -8&\hskip1.3ex0\\\end{psmallmatrix}+
\begin{psmallmatrix}-8&-8\\ \hskip1.3ex 8&\hskip1.3ex 8\\\end{psmallmatrix}+2\begin{psmallmatrix}4&2\\ 0&0\\\end{psmallmatrix}
\right\}\\
&=8,
\end{aligned}
$$
and accounting for the orientation and $-2$ in the expression above for the WP form yields the asserted constant
$\gamma=-4\omega$.

\medskip 


\noindent Finally in the remaining case to compute $\gamma(\bar\vartheta_A,\bar\vartheta_{TA})$, replace $A\in {\rm PSL}_2$ by {$SU\hskip-.55exA$}, so that
$$\gamma(\bar\vartheta_{SU\hskip-.55exA},\bar\vartheta_{TSU\hskip-.55exA})=\gamma(\bar\vartheta_{U\hskip-.55exA},\bar\vartheta_A)=+8$$
by skew symmetry since $\bar\vartheta _{SB}=\bar\vartheta_B$, for $B\in {\rm PSL}_2$, and
$$TSU\hskip-.55exA=TT^{-1}UT^{-1}U\hskip-.55exA=UT^{-1}UA=S\hskip-.35exA.$$\end{proof}

\bigskip

\section{Universal automorphic 1-form on ${\rm PPSL}_2({\mathbb R})$}

\bigskip

\noindent Given a decoration $\tilde\tau$ on a tesselation
$\tau$ with framing
framing $\mathcal F$, define
$$\begin{aligned}
\xi_{\mathcal F}&
=
{1\over 2}\sum_{A\in{\rm PSL}_2} \bar\vartheta_A~d{\rm log}\,\lambda _{e_A}
\in\Omega^1(\widetilde{{{\rm Homeo}}_+^n},ppsl_2),
\end{aligned}$$
where $\bar\vartheta_A$ is the wavelet $\vartheta_A$ normalized with respect to $\mathcal F$, and $\lambda_e$ denotes the lambda length of $e=e_A$.  Alternatively, we may write $\xi_{\mathcal F}=\sum_{e\in\tau_*} \bar\vartheta_e~d{\rm log}\,\lambda _e$ since $\vartheta_A=\vartheta_{SA}$ and $\lambda_e=\lambda_{e_A}=\lambda_{e_{SA}}$.

\medskip

\noindent $\xi_{\mathcal F}$ is to be apprehended as a 1-form on the group $\widetilde{{\rm Homeo}^n_+}$
taking values in its Lie algebra $ppsl_2$.  It evidently could be interpreted as
the Maurer-Cartan form since the result of applying it to a multiplicative deformation of the lambda length, or in other words an additive deformation of the logarithm of the lambda length, is precisely the corresponding normalized vector field on ${\mathbb S}^1$ in its Lie algebra $ppsl_2$, as was discussed in Section \ref{sec:tess}.

\medskip

\noindent More precisely, we have the diagram
$$\begin{aligned}
&\widetilde{{\rm Homeo}_+^n}\approx\widetilde{{\rm Homeo}_+^n}/{M\ddot ob}\\
&\hskip 3ex\downarrow\pi\hskip 10ex \downarrow\\
{\rm Homeo}_+\stackrel[s_{\mathcal F}]{p}{\rightleftarrows}~&{\rm Homeo}_+^n\approx {\rm Homeo}_+/{M\ddot ob},
\end{aligned}$$
where $\pi$ is the forgetful map, $p$ is the projection
given by precomposition with
$L_{f({0\over 1}),f({1\over 0}),f({1\over 1})}^{{0\over 1},{1\over 0},{1\over 1}}$,
and for any framing ${\mathcal F}=(u,v,w)$, the section $s_{\mathcal F}$ of $p$ is given by precomposition with $L_{{0\over 1},{1\over 0},{1\over 1}}^{u,v,w}.$
It is a tautology that if
$\omega$ is the Maurer-Cartan form of ${\rm Homeo}_+$, then
 $\xi_{\mathcal F}=\pi^*\circ s^*_{\mathcal F}\, \omega$.

\begin{theorem}\label{thm:MCflip} For any fixed framing $\mathcal F\in C_3$, 
the Lie-algebra valued 1-form $\xi_{\mathcal F}\in\Omega^1(\widetilde{{\rm Homeo}}_+,ppsl_2)$ is invariant under flips.
\end{theorem}

\begin{proof} As a Lie algebra valued 1-form on a group, it is described by right-translating from the cotangent plane of the identity, and our computations
will take place there.

\medskip

\noindent
Recall that the Farey tesselation $\tau_*$ admits a canonical decoration $\tilde\tau_*$ determined by the property that all the lambda lengths are constant equal to unity, and that this is the identity element of the group $\widetilde{{\rm Homeo}_+^n}$ as discussed in Section \ref{sec:tess}. 

\medskip

\noindent To establish notation, refer to Figure 1 and consider the four generation-zero and -one Farey points
${1\over 0},-{1\over 1},{0\over 1}, {1\over 1}$  decomposing ${\mathbb S}^1$ into four circular intervals which are conflated with the respective quadrants I,II,III and IV containing them.  The convex hull of these points is an ideal quadrilateral with oriented frontier edges $a=e_{ST},b=e_{SU},c=e_{ U^{-1}}$ and $d=e_{ T^{-1}}$ in this counter-clockwise order
starting from ${1\over 0}$.  The doe $e_I$ will be denoted simply $e=e_I$, and the edge arising from its flip, with respective initial and terminal points
${1\over 1}$ and $-{1\over 1}$, will be denoted $f$.

\medskip 

\noindent The second-generation 
$-{2\over 1},-{1\over 2},{1\over 2},{2\over 1}$ Farey points further decompose the circular intervals ${\rm I}={\rm I}_-\cup {\rm I}_+,\ldots ,{\rm IV}={\rm IV}_-\cup {\rm IV}_+$ occurring in this counter-clockwise order starting from ${1\over 0}$.  The Farey points of generation at most two together span an ideal octagon with oriented frontier edges $e_{ST^2}$,  $e_{SUT}$,  $e_{TSU}$, $e_{SU^2}$, $e_{ U^{-2}}$, $e_{ T^{-1} U^{-1}}$, $e_{ U^{-1} T^{-1}}$, and $e_{ T^{-2}}$ in this counter-clockwise order starting from
${1\over 0}$, as indicated in Figure 2.

\medskip

\noindent The hyperbolic transformations which are primitive in PSL$_2$ along the respective axes
$e$ and $f$ are given by the exponentials of $\begin{psmallmatrix}-1& 0\\\hskip1.3ex0&1\\\end{psmallmatrix}$ and $\begin{psmallmatrix}0& 1\\1&0\\\end{psmallmatrix}$, and along
each of the other edges $e_A$ above by the exponential of $X_A\in sl_2$ with
$$X_{ST}=\begin{psmallmatrix}1& \hskip1.4ex 2\\0&-1\\\end{psmallmatrix},~~X_{SU}=\begin{psmallmatrix}\hskip 1.4ex1& \hskip1.4ex 0\\-2&-1\\\end{psmallmatrix},~~X_{ U^{-1}}=\begin{psmallmatrix}-1& \hskip 1.4ex 0\\-2&\hskip1.3ex 1\\\end{psmallmatrix},~~X_{ T^{-1}}=\begin{psmallmatrix}-1& 2\\\hskip 1.3ex 0&1\\\end{psmallmatrix},$$
and
$$\begin{aligned}X_{ST^2}&=\begin{psmallmatrix}1& \hskip 1.2ex 4\\0&-1\\\end{psmallmatrix},~~
X_{SUT}=\begin{psmallmatrix}\hskip1.3ex 3&\hskip 1.3ex 4\\-2&-3\\\end{psmallmatrix},\\~~
X_{TSU}&=\begin{psmallmatrix}\hskip 1.3ex 3&\hskip1.2ex 2\\-4&-3\\\end{psmallmatrix},~~
X_{SU^2}=\begin{psmallmatrix}\hskip1.4ex1& \hskip1.4ex 0\\-4&-1\\\end{psmallmatrix},
\\
X_{ U^{-2}}&=\begin{psmallmatrix}-1&  0\\-4&1\\\end{psmallmatrix},~~
X_{ T^{-1} U^{-1}}=\begin{psmallmatrix}-3& 2\\-4&3\\\end{psmallmatrix},\\
X_{ U^{-1} T^{-1}}&=\begin{psmallmatrix}-3& 4\\-2&3\\\end{psmallmatrix},~~
X_{ U^{-1} T^{-1}}=\begin{psmallmatrix}-1& 4\\\hskip 1.2ex 0&1\\\end{psmallmatrix}.
\end{aligned}$$
The proof is a direct computation and is omitted. 

\medskip

\noindent 
These are the logarithms of the component earthquakes for the wavelets
$\vartheta_x$, for $x\in\{a,b,c,d,e,f\}$, so it is not difficult to combine them
four at a time and normalize with respect to the the standard framing
${\mathcal F}=({0\over 1},{1\over 0},{1\over 1})$ on $\tau _*$.
One finds
$$\begin{aligned}
\bar\vartheta_e=\vartheta_e&=\begin{cases}
\begin{psmallmatrix}-1&-2\\\hskip 1.2ex 0&\hskip 1.3ex 1\\ \end{psmallmatrix},&{\rm on~I};\\
\begin{psmallmatrix}\hskip 1.3ex1&\hskip 1.2ex 0\\-2&-1\\ \end{psmallmatrix},&{\rm on~II};\\
\hskip 1.4ex\begin{psmallmatrix}1&\hskip1.2ex0\\2&-1\\ \end{psmallmatrix},&{\rm on~III};\\
\hskip 1.4ex\begin{psmallmatrix}-1&2\\\hskip 1.2ex0&1\\ \end{psmallmatrix},&{\rm on~IV},\\
\end{cases}
\end{aligned}$$
$$\begin{aligned}
\vartheta_a&=\begin{cases}
\begin{psmallmatrix}-1&-4\\\hskip1.4ex0&\hskip1.4ex1\\ \end{psmallmatrix},&{\rm on~I_-};\\
\begin{psmallmatrix}\hskip1.4ex 3&\hskip.9ex 4\\-2&-3\\ \end{psmallmatrix},&{\rm on~I_+};\\
\hskip 1.1ex\begin{psmallmatrix}-1&0\\-2&1\\ \end{psmallmatrix},&{\rm on~II};\\
\hskip 1.2ex\begin{psmallmatrix}-1&0\\\hskip1.4ex0&1\\ \end{psmallmatrix},&{\rm on~III~and~IV},\\
\end{cases}
~~~~{\rm so}~~~~
\bar\vartheta_a=\begin{cases}
\hskip1.6ex\begin{psmallmatrix}0&-4\\0&\hskip1.2ex0\\ \end{psmallmatrix},&{\rm on~I_-};\\
\begin{psmallmatrix}\hskip1.4ex4&\hskip1.4ex4\\-2&-4\\ \end{psmallmatrix},&{\rm on~I_+};\\
\hskip 2.5ex\begin{psmallmatrix}0&0\\2&0\\ \end{psmallmatrix},&{\rm on~II};\\
\hskip 2.6ex\begin{psmallmatrix}0&0\\0&0\\ \end{psmallmatrix},&{\rm on~III~and~IV},\\
\end{cases}
\end{aligned}$$

$$\begin{aligned}
\vartheta_b&=\begin{cases}
\hskip1.1ex\begin{psmallmatrix}1&\hskip1.4ex2\\0&-1\\ \end{psmallmatrix},&{\rm on~I};\\
\begin{psmallmatrix}-3&-2\\\hskip1.2ex 4&\hskip 1.4ex 3\\ \end{psmallmatrix},&{\rm on~II_-};\\
\begin{psmallmatrix}\hskip1.4ex 1&\hskip1.3ex 0\\-4&-1\\ \end{psmallmatrix},&{\rm on~II_+};\\
\hskip 1.3ex\begin{psmallmatrix}1&\hskip1.3ex 0\\0&-1\\ \end{psmallmatrix},&{\rm on~III~and~IV},\\
\end{cases}
~~~~{\rm so}~~~~
\bar\vartheta_b=\begin{cases}
\hskip2.4ex\begin{psmallmatrix}0&2\\0&0\\ \end{psmallmatrix},&{\rm on~I};\\
\begin{psmallmatrix}-4&-2\\\hskip1.2ex 4&\hskip 1.4ex 4\\ \end{psmallmatrix},&{\rm on~II_-};\\
\hskip 1.25ex\begin{psmallmatrix}\hskip1.4ex 0&0\\-4&0\\ \end{psmallmatrix},&{\rm on~II_+};\\
\hskip 2.4ex\begin{psmallmatrix}0&0\\0&0\\ \end{psmallmatrix},&{\rm on~III~and~IV},\\
\end{cases}
\end{aligned}$$
$$\begin{aligned}
\bar\vartheta_c=\vartheta_c&=\begin{cases}
\begin{psmallmatrix}1&\hskip1.4ex 0\\0&-1\\ \end{psmallmatrix},&{\rm on~I~and~II};\\
\begin{psmallmatrix}1&\hskip1.4ex 0\\4&-1\\ \end{psmallmatrix},&{\rm on~III_-};\\
\begin{psmallmatrix}-3&2\\-4&3\\ \end{psmallmatrix},&{\rm on~III_+};\\
\begin{psmallmatrix}1&-2\\0&-1\\ \end{psmallmatrix},&{\rm on~IV},\\
\end{cases}~~{\rm and}~~
\bar\vartheta_d=\vartheta_d=\begin{cases}
\begin{psmallmatrix}-1&0\\\hskip1.4ex 0&1\\ \end{psmallmatrix},&{\rm on~I~and ~II};\\
\begin{psmallmatrix}-1&0\\-2&1\\ \end{psmallmatrix},&{\rm on~III};\\
\begin{psmallmatrix}3&-4\\2&-3\\ \end{psmallmatrix},&{\rm on~IV_-};\\
\begin{psmallmatrix}-1&4\\\hskip 1.4ex0&1\\ \end{psmallmatrix},&{\rm on~IV_+}\\
\end{cases}
\end{aligned}$$

\noindent As illustrated in Figure 5, letting $\tilde x=d{\rm log} \,\lambda_x$, for $x\in\{a,b,c,d,e,f\}$, and putting all of this together, we find that
$$\begin{aligned}
&\sum_{x\in{\{a,b,c,d,e\}}}\bar\vartheta_x~d{\rm log}\,\lambda_x\\
&\hskip1cm=~~\begin{cases}

\smallskip
\begin{psmallmatrix}0&-4\\0&\hskip1.2ex 0\end{psmallmatrix}\tilde a
+\begin{psmallmatrix}0&2\\0&0\end{psmallmatrix}\tilde b
+\begin{psmallmatrix}1&\hskip1.2ex0\\0&-1\end{psmallmatrix}\tilde c
+\begin{psmallmatrix}-1&0\\\hskip1.2ex 0&1\end{psmallmatrix}\tilde d
+\begin{psmallmatrix}-1&-2\\\hskip1.2ex 0&\hskip1.4ex 1\end{psmallmatrix}\tilde e,&{\rm on}~{\rm I}_-;\\

\smallskip

\begin{psmallmatrix}\hskip1.2ex4&\hskip1.2ex4\\-2&-4\end{psmallmatrix}\tilde a
+\begin{psmallmatrix}0&2\\0&0\end{psmallmatrix}\tilde b
+\begin{psmallmatrix}1&\hskip1.2ex0\\0&-1\end{psmallmatrix}\tilde c
+\begin{psmallmatrix}-1&0\\\hskip1.2ex 0&1\end{psmallmatrix}\tilde d
+\begin{psmallmatrix}-1&-2\\\hskip1.2ex 0&\hskip1.4ex 1\end{psmallmatrix}\tilde e,&{\rm on}~{\rm I}_+;\\

\smallskip

\begin{psmallmatrix}0&0\\2&0\end{psmallmatrix}\tilde a
+\begin{psmallmatrix}-4&-2\\\hskip1.2ex4&\hskip1.2ex4\end{psmallmatrix}\tilde b
+\begin{psmallmatrix}1&\hskip1.2ex 0\\0&-1\end{psmallmatrix}\tilde c
+\begin{psmallmatrix}-1&0\\\hskip1.2ex 0&1\end{psmallmatrix}\tilde d
+\begin{psmallmatrix}\hskip1.4ex1&\hskip1.2ex 0\\-2&-1\end{psmallmatrix}\tilde e,&{\rm on}~{\rm II}_-;\\

\smallskip

\begin{psmallmatrix}0&0\\2&0\end{psmallmatrix}\tilde a
+\begin{psmallmatrix}\hskip1.2ex0&0\\-4&0\end{psmallmatrix}\tilde b
+\begin{psmallmatrix}1&\hskip1.2ex 0\\0&-1\end{psmallmatrix}\tilde c
+\begin{psmallmatrix}-1&0\\\hskip1.2ex 0&1\end{psmallmatrix}\tilde d
+\begin{psmallmatrix}\hskip1.4ex1&\hskip1.2ex 0\\-2&-1\end{psmallmatrix}\tilde e,&{\rm on}~{\rm II}_+;\\
\smallskip

\begin{psmallmatrix}1&\hskip1.2ex0\\4&-1\end{psmallmatrix}\tilde c
+\begin{psmallmatrix}-1&0\\-2&1\end{psmallmatrix}\tilde d
+\begin{psmallmatrix}1&\hskip1.2ex0\\2&-1\end{psmallmatrix}\tilde e,&{\rm on}~{\rm III}_-;\\

\smallskip

\begin{psmallmatrix}-3&2\\-4&3\end{psmallmatrix}\tilde c
+\begin{psmallmatrix}-1&0\\-2&1\end{psmallmatrix}\tilde d
+\begin{psmallmatrix}1&\hskip1.2ex0\\2&-1\end{psmallmatrix}\tilde e,&{\rm on}~{\rm III}_+;\\

\smallskip

\begin{psmallmatrix}1&-2\\0&-1\end{psmallmatrix}\tilde c
+\begin{psmallmatrix}3&-4\\2&-3\end{psmallmatrix}\tilde d
+\begin{psmallmatrix}-1&2\\\hskip1.2ex 0&1\end{psmallmatrix}\tilde e,&{\rm on}~{\rm IV}_-;\\

\smallskip

\begin{psmallmatrix}1&-2\\0&-1\end{psmallmatrix}\tilde c
+\begin{psmallmatrix}-1&4\\\hskip1.2ex0&1\end{psmallmatrix}\tilde d
+\begin{psmallmatrix}-1&2\\\hskip1.2ex 0&1\end{psmallmatrix}\tilde e,&{\rm on}~{\rm IV}_+.\\

\end{cases}
\end{aligned}$$

\setcounter{figure}{4}
\captionsetup[figure]{font=small,skip=0pt}
\begin{center}
\includegraphics[scale=.6]{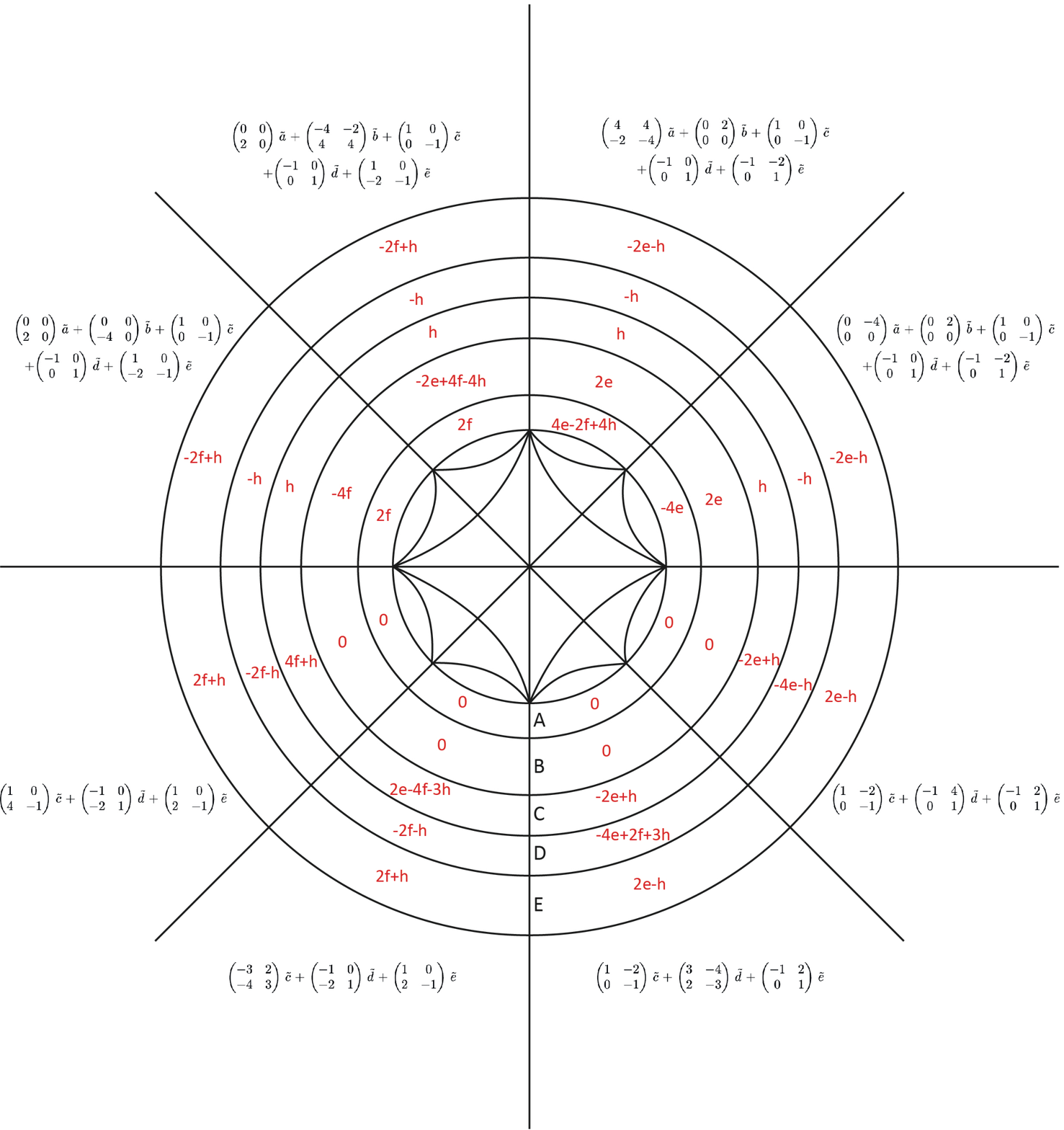}
\captionof{figure}{Depiction of $\bar\vartheta_a,\ldots, \bar\vartheta_e$ in the respective rings
$A,\ldots,E$ and $\sum_{x\in{\{a,b,c,d,e\}}}\bar\vartheta_x~d{\rm log}\,\lambda_x$ on the outside.}
\end{center}

\noindent Let $\tau'$ denote the tesselation that arises from $\tau_*$ by a flip along the edge $e\in\tau_*$, let $\vartheta_x'$ denote the vectorfield on ${\mathbb S}^1$
corresponding in the natural way to the edge $x\in\{ a,b,c,d,f\}\subseteq\tau'$,
and let $\bar\vartheta_x'$ denote the normalization of $\vartheta_x'$ relative to the
same framing $\mathcal F$.  Thus, we have
$$\bar\vartheta_f'=\vartheta'_f=\bar\vartheta_e,$$ and a further consequence of
the primitives for logarithms of hyperbolic transformations computed before is that

\captionsetup[figure]{font=small,skip=-0pt}
\begin{center}
\includegraphics[scale=.6]{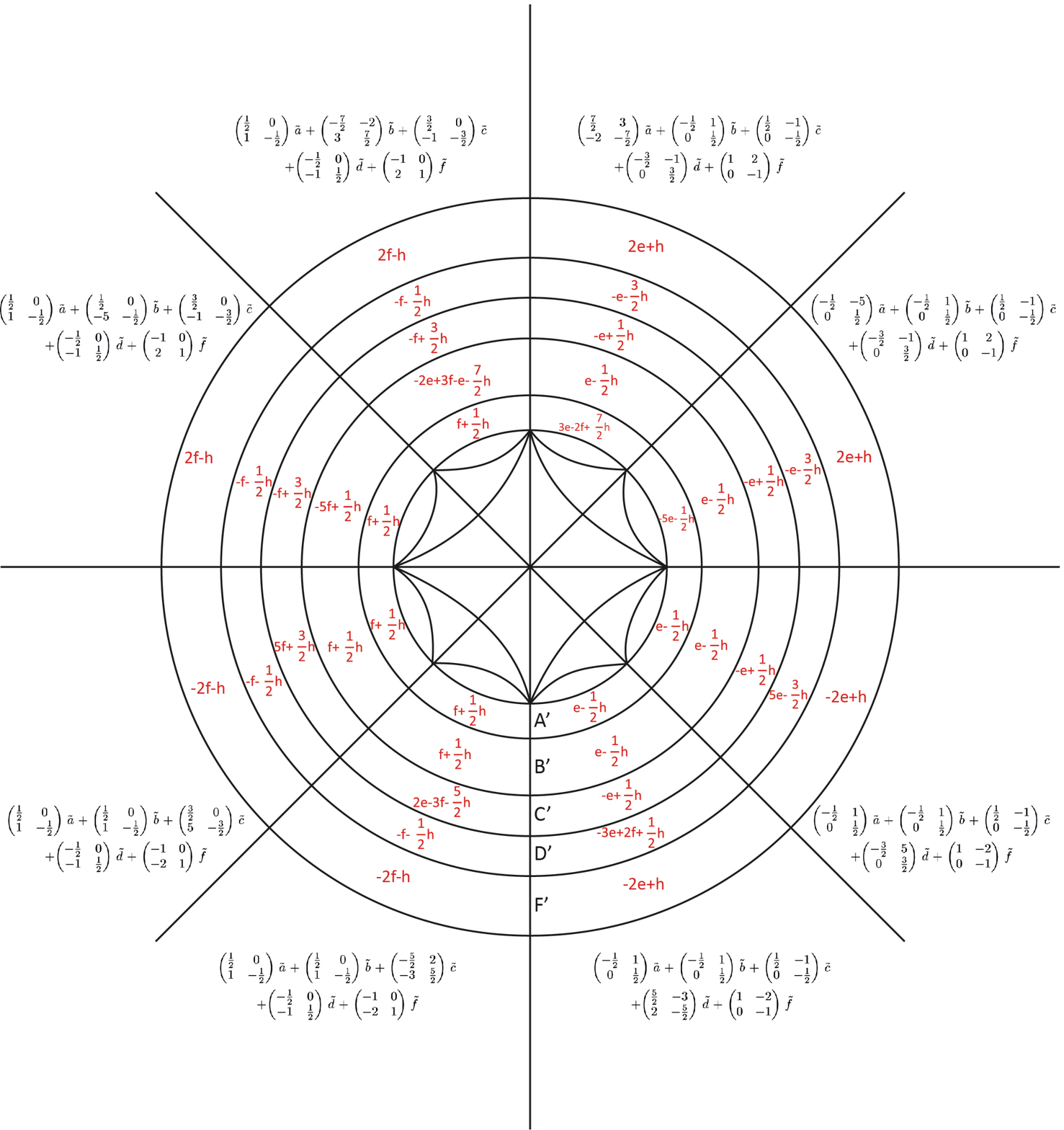}
\captionof{figure}{Depiction of $\bar\vartheta_a',\ldots, \bar\vartheta_f'$ in the respective rings
$A',\ldots D', F'$ and $\sum_{x\in{\{a,b,c,d,f\}}}\bar\vartheta_x'~d{\rm log}\,\lambda_x$ on the outside.}
\end{center}

$$\begin{aligned}
\bar\vartheta_a'=\vartheta_a'+\begin{psmallmatrix}
{1\over 2}&-1\\0&-{1\over 2}\end{psmallmatrix},~~{\rm where}~~~\vartheta_a'&=\begin{cases}
\begin{psmallmatrix}-1&-4\\\hskip1.2ex 0&\hskip 1.4ex1\\ \end{psmallmatrix},&{\rm on~I_-};\\
\begin{psmallmatrix}\hskip1.2ex 3&\hskip1.2ex4\\-2&-3\\ \end{psmallmatrix},&{\rm on~I_+};\\
\hskip 2.5ex\begin{psmallmatrix}0&1\\1&0\\ \end{psmallmatrix},&{\rm on~II~and~III};\\
\hskip 1.2ex\begin{psmallmatrix}-1&2\\\hskip1.4ex0&1\\ \end{psmallmatrix},&{\rm on~IV},
\end{cases}
\end{aligned}$$
$$\begin{aligned}
\bar\vartheta_b'=\vartheta_b'+
\begin{psmallmatrix}
-{1\over 2}&0\\-1&{1\over 2}\end{psmallmatrix},~~{\rm where}~~~\vartheta_b'&=\begin{cases}
\hskip 2.3ex\begin{psmallmatrix}0&1\\1&0\\ \end{psmallmatrix},&{\rm on~I~and~IV};\\
\begin{psmallmatrix}-3&-2\\\hskip1.2ex4&\hskip1.2ex3\\ \end{psmallmatrix},&{\rm on~II_-};\\
\begin{psmallmatrix}\hskip1.4ex1&\hskip1.0ex0\\-4&-1\\ \end{psmallmatrix},&{\rm on~II_+};\\
\hskip1.2ex\begin{psmallmatrix}1&\hskip 1.2ex 0\\2&-1\\ \end{psmallmatrix},&{\rm on~III},\\
\end{cases}
\end{aligned}$$
$$\begin{aligned}
\bar\vartheta_c'=\vartheta_c'+
\begin{psmallmatrix}
{1\over 2}&\hskip1.2ex0\\1&-{1\over 2}\end{psmallmatrix},~~{\rm where}~~~\vartheta_c'&=\begin{cases}
\begin{psmallmatrix}\hskip1.2ex0&-1\\-1&\hskip1.2ex0\\ \end{psmallmatrix},&{\rm on~I~and~IV};\\
\begin{psmallmatrix}\hskip 1.4ex 1&\hskip1.2ex0\\-2&-1\\ \end{psmallmatrix},&{\rm on~II};\\
\hskip1.4ex\begin{psmallmatrix}1&\hskip1.2ex0\\4&-1\\ \end{psmallmatrix},&{\rm on~III_-};\\
\begin{psmallmatrix}-3&\hskip1.2ex2\\-4&-1\\ \end{psmallmatrix},&{\rm on~III_+},\\
\end{cases}
\end{aligned}$$
$$\begin{aligned}
\bar\vartheta_d'=\vartheta_d'+
\begin{psmallmatrix}
-{1\over 2}&1\\\hskip1.2ex 0&{1\over 2}\end{psmallmatrix},~~{\rm where}~~~\vartheta_d'&=\begin{cases}
\begin{psmallmatrix}-1&-2\\\hskip1.2ex0&\hskip1.4ex1\\ \end{psmallmatrix},&{\rm on~I,};\\
\begin{psmallmatrix}\hskip 1.2ex 0&-1\\-1&\hskip1.2ex 0\\ \end{psmallmatrix},&{\rm on~II~and~III};\\
\hskip1.2ex\begin{psmallmatrix}-3&4\\-2&3\\ \end{psmallmatrix},&{\rm on~IV_-};\\
\hskip1.4ex\begin{psmallmatrix}-{1}&4\\\hskip1.2ex 0&1\\ \end{psmallmatrix},&{\rm on~IV_+}.\\
\end{cases}
\\\\
\end{aligned}$$

\noindent Again putting all this together as illustrated in Figure 6, we find
$$\begin{aligned}
&\sum_{x\in{\{a,b,c,d,f\}}}\bar\vartheta_x'~d{\rm log}\,\lambda_x\\
&\hskip1cm=~~\begin{cases}

\medskip
\begin{psmallmatrix}-{1\over 2}&-5\\\hskip 1.2ex0&\hskip1.2ex {1\over 2}\end{psmallmatrix}\tilde a
+\begin{psmallmatrix}-{1\over 2}&1\\\hskip1.2ex0&{1\over 2}\end{psmallmatrix}\tilde b
+\begin{psmallmatrix}{1\over 2}&-1\\0&-{1\over 2}\end{psmallmatrix}\tilde c
+\begin{psmallmatrix}-{3\over 2}&-1\\\hskip1.2ex 0&\hskip 1.2ex{3\over 2}\end{psmallmatrix}\tilde d
+\begin{psmallmatrix}1&\hskip 1.2ex2\\0&-1\end{psmallmatrix}\tilde f,&{\rm on}~{\rm I}_-;\\

\medskip

\begin{psmallmatrix}\hskip 1.3ex{7\over 2}&\hskip 1.4ex3\\-2&-{1\over 2}\end{psmallmatrix}\tilde a
+\begin{psmallmatrix}-{1\over 2}&1\\\hskip 1.2ex 0&{1\over 2}\end{psmallmatrix}\tilde b
+\begin{psmallmatrix}{1\over 2}&-1\\0&-{1\over 2}\end{psmallmatrix}\tilde c
+\begin{psmallmatrix}-{3\over 2}&-1\\\hskip1.2ex 0&\hskip 1.2ex{3\over 2}\end{psmallmatrix}\tilde d
+\begin{psmallmatrix}1&\hskip1.2ex 2\\0&-1\end{psmallmatrix}\tilde f,&{\rm on}~{\rm I}_+;\\

\medskip

\begin{psmallmatrix}{1\over 2}&\hskip 1.2ex0\\1&-{1\over 2}\end{psmallmatrix}\tilde a
+\begin{psmallmatrix}-{7\over 2}&-2\\\hskip1.2ex3&\hskip1.2ex{7\over 2}\end{psmallmatrix}\tilde b
+\begin{psmallmatrix}\hskip 1.2ex{3\over 2}&\hskip1.2ex 0\\-1&-{3\over 2}\end{psmallmatrix}\tilde c
+\begin{psmallmatrix}-{1\over 2}&0\\-1&{1\over 2}\end{psmallmatrix}\tilde d
+\begin{psmallmatrix}-1&0\\\hskip 1,2ex 2&1\end{psmallmatrix}\tilde f,&{\rm on}~{\rm II}_-;\\

\medskip

\begin{psmallmatrix}{1\over 2}&\hskip 1.2ex0\\1&-{1\over 2}\end{psmallmatrix}\tilde a
+\begin{psmallmatrix}\hskip1.2ex{1\over 2}&\hskip 1.2ex0\\-5&-{1\over 2}\end{psmallmatrix}\tilde b
+\begin{psmallmatrix}\hskip 1.2ex{3\over 2}&\hskip1.4ex 0\\-1&-{3\over 2}\end{psmallmatrix}\tilde c
+\begin{psmallmatrix}-{1\over 2}&0\\-1&{1\over 2}\end{psmallmatrix}\tilde d
+\begin{psmallmatrix}-1&0\\\hskip 1.2ex 2&1\end{psmallmatrix}\tilde f,&{\rm on}~{\rm II}_+;\\

\medskip

\begin{psmallmatrix}{1\over 2}&\hskip 1.2ex 0\\1&-{1\over 2}\end{psmallmatrix}\tilde a
+\begin{psmallmatrix}{1\over 2}&\hskip1.2ex 0\\1&-{1\over 2}\end{psmallmatrix}\tilde b
+\begin{psmallmatrix}{3\over 2}&\hskip 1.2ex 0\\5&-{3\over 2}\end{psmallmatrix}\tilde c
+\begin{psmallmatrix}-{1\over 2}&0\\-1&{1\over 2}\end{psmallmatrix}\tilde d
+\begin{psmallmatrix}-1&0\\-2&1\end{psmallmatrix}\tilde f,&{\rm on}~{\rm III}_-;\\

\medskip

\begin{psmallmatrix}{1\over 2}&\hskip 1.2ex 0\\1&-{1\over 2}\end{psmallmatrix}\tilde a
+\begin{psmallmatrix}{1\over 2}&\hskip 1.2ex0\\1&-{1\over 2}\end{psmallmatrix}\tilde b
+\begin{psmallmatrix}-{5\over 2}&2\\-3&{5\over 2}\end{psmallmatrix}\tilde c
+\begin{psmallmatrix}-{1\over 2}&0\\-1&{1\over 2}\end{psmallmatrix}\tilde d
+\begin{psmallmatrix}-1&\hskip 1.2ex0\\-2&-1\end{psmallmatrix}\tilde f,&{\rm on}~{\rm III}_+;\\

\medskip

\begin{psmallmatrix}-{1\over 2}&1\\\hskip 1.2ex 0&{1\over 2}\end{psmallmatrix}\tilde a
+\begin{psmallmatrix}-{1\over 2}&1\\\hskip1.2ex0&{1\over 2}\end{psmallmatrix}\tilde b
+\begin{psmallmatrix}{1\over 2}&-1\\0&-{1\over 2}\end{psmallmatrix}\tilde c
+\begin{psmallmatrix}{5\over 2}&-3\\2&-{5\over 2}\end{psmallmatrix}\tilde d
+\begin{psmallmatrix}1&-2\\0&-1\end{psmallmatrix}\tilde f,&{\rm on}~{\rm IV}_-;\\

\medskip

\begin{psmallmatrix}-{1\over 2}&1\\\hskip 1.2ex0&{1\over 2}\end{psmallmatrix}\tilde a
+\begin{psmallmatrix}-{1\over 2}&1\\\hskip 1.2ex 0&{1\over 2}\end{psmallmatrix}\tilde b
+\begin{psmallmatrix}{1\over 2}&-1\\0&-{1\over 2}\end{psmallmatrix}\tilde c
+\begin{psmallmatrix}-{3\over 2}&5\\\hskip1.2ex0&{3\over 2}\end{psmallmatrix}\tilde d
+\begin{psmallmatrix}1&-2\\0&-1\end{psmallmatrix}\tilde f,&{\rm on}~{\rm IV}_+;
\smallskip

\end{cases}
\end{aligned}
$$

\noindent Using the fact that
$$\tilde e+\tilde f={1\over{ac+bd}}(\tilde a+\tilde b+\tilde c+\tilde d)=
{1\over 2}(\tilde a+\tilde b+\tilde c+\tilde d)$$
by the Ptolemy relation $ef=ac+bd$, since
$a=b=c=d=1$ for the identity element $\tilde\tau_*$,
a computation with rather miraculous cancellations as illustrated in Figure 7
confirms that
$$\sum_{x\in{\{a,b,c,d,e\}}}\bar\vartheta_x~d{\rm log}\,\lambda_x
~~=~~\sum_{x\in{\{a,b,c,d,e\}}}\bar\vartheta_x'~d{\rm log}\,\lambda_x,$$
as required to check invariance under the flip along the doe for the standard framing ${\mathcal F}=({0\over 1},{1\over 0},{1\over 1})$.


\captionsetup[figure]{font=small,skip=-0in}
\begin{center}
\includegraphics[trim=0 0 0 20,clip,width=.85\textwidth]{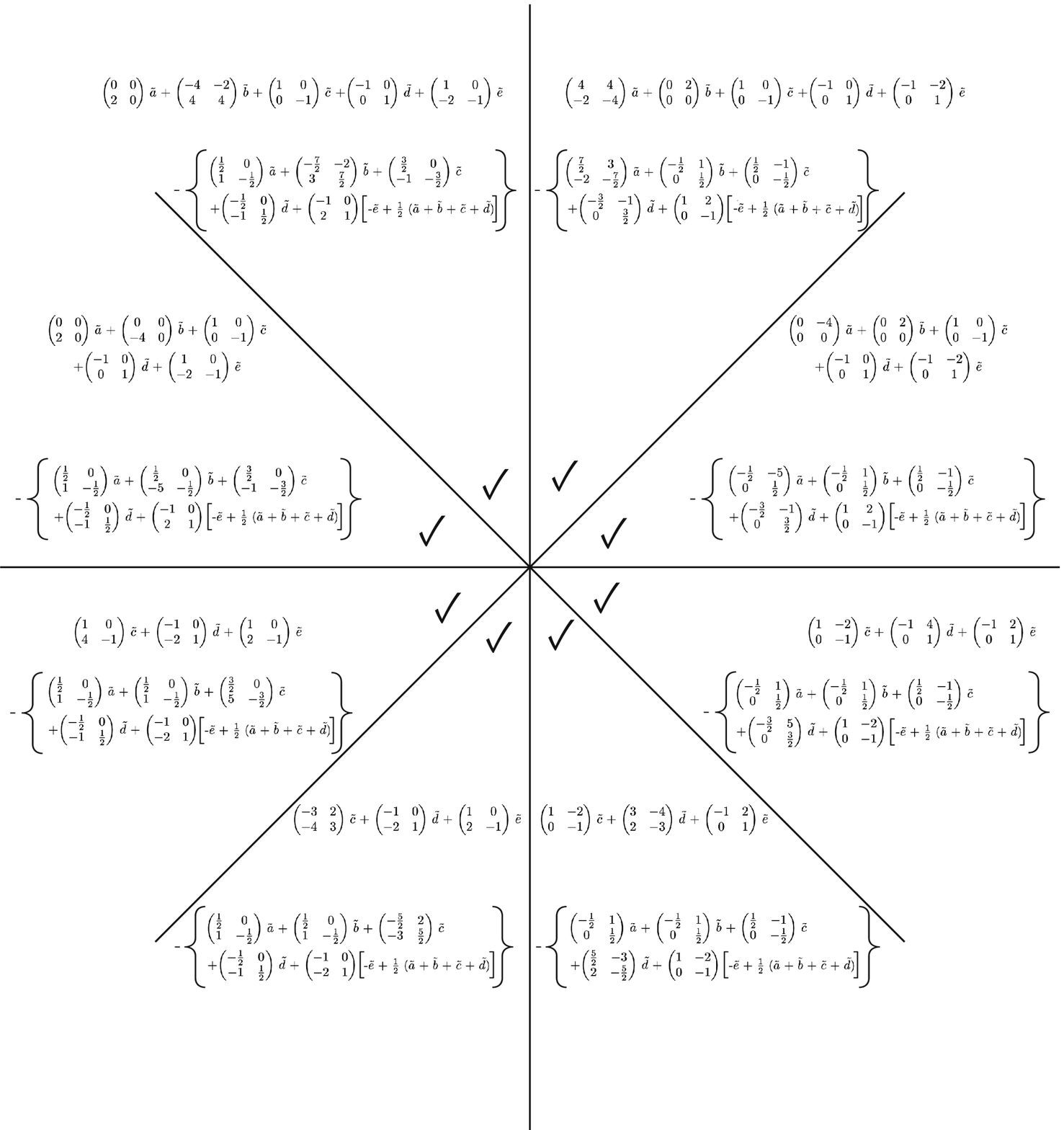}
\captionof{figure}{Difference between expressions in Figures 5 and 6 vanishes
for all lambda lengths equal to unity, i.e., at the universal punctured arithmetic surface.}
\end{center}

\noindent This establishes the invariance for a flip on the doe with the associated framing.
In the general case first of all, the 1-form is invariant under push-forward by the M\"obius group since lambda lengths are $M\ddot ob$-invariant, so we may assume that the framing is one that is associated to a doe.  Likewise without loss of generality, we may assume that the flip is performed on the edge with endpoints ${0\over 1},{1\over 0}\in{\mathbb S}^1$, with the doe however potentially located elsewhere.  

\medskip

\noindent Thus, the asserted invariance must be checked for the various possibilities
of relative position of doe and the edge with endpoints ${0\over 1},{1\over 0}$ to be flipped.  In the presented computation, these two edges coincide, and there are four further cases 
to check depending upon which quadrant in ${\mathbb C}\supseteq{\mathbb D}$
contains the doe, thereby determining the normalization.
These four similar but simpler calculations, each corresponding to a different normalization of the common unnormalized formulae already presented in the proof, are left as an exercise for the reader.
\end{proof}

\bigskip

\noindent
Consider an arbitrary framing ${\mathcal F}$ and let 
$$\eta_{\mathcal F}=s^*_{\mathcal F}\,\omega\in\Omega^1({\rm Homeo}_+/{M\ddot ob}),$$
which by general principles should satisfy
$$\begin{aligned}
&\underline{{\rm Maurer-Cartan~equation}}:~0=d\eta_{\mathcal F}+{{1\over 2}}[\eta_{\mathcal F},\eta_{\mathcal F}],~{\rm for~any~framing}~{\mathcal F};\\
\smallskip 
&\underline{{\rm Compatibility~Equation}}: \eta_{\mathcal G}={\rm Ad}(h^{-1}_{{\mathcal F}{\mathcal G}})\,\eta_{\mathcal F}~+~h^*_{{\mathcal F}{\mathcal G}}\,\omega_{M\ddot ob},
\mbox{for framings }\\ 
&{\mathcal F,G},
{\rm where }~h_{{\mathcal F}{\mathcal G}}=s_{\mathcal G}\circ s_{\mathcal F}^{-1}, \mbox{Ad denotes the adjoint} {\rm ~and}~\omega_{M\ddot ob}~{\rm is~the}\\
& \mbox{Maurer-Cartan form on}~{M\ddot ob}. \\
\end{aligned}$$

\medskip

\noindent Note that the Compatibility Equation together with the result for the flip on the doe with its associated framing would give an alternative proof of Theorem \ref{thm:MCflip}.
This putative proof is a bit of a swindle because {\sl there  is
no Maurer-Cartan form on the topological group} ${\rm Homeo}_+$ except for the version here, and hence no a priori Compatibility Equation.

\begin{remark} The expression for $\xi_{\mathcal F}$ as a Poincar\'e series in framed holographic coordinates is easily derived from Remark \ref{rem:stablam}.
\end{remark}

\medskip

\noindent As we have noted, the form $\xi_{\mathcal F}$ is morally the Maurer-Cartan form in the sense that its value on the tangent vector ${{\partial}\over{\partial {\rm ln}\lambda_A}}$ at the identity of ${\rm PPSL}_2({\mathbb R})$ is the corresponding normalized element $\bar\vartheta_A$ of the Lie algebra $ppsl_2$.  However, $\xi_{\mathcal F}$ has been shown
here to be invariant under the lattice
${\rm PPSL}_2({\mathbb Z})$ rather than
under the full group ${\rm PPSL}_2({\mathbb R})$, whose invariance we may conjecture.
On the other hand, one might therefore more squarely regard $\xi_{\mathcal F}$ solely as an automorphic form, and this is the viewpoint of the next section.

\section{Future perspectives}

\medskip

\noindent
We conclude with a heuristic discussion of constructing the automorphic representation for our universal triad corresponding to the automorphic 1-form obtained
in Section 7.  To this end, we combine the original constructions of \cite{VGG} and of
Vaughan Jones \cite{jones1,jones2}, also described in \cite{OS}.  

\medskip

\noindent
We recall from Appendix B that we have constructed an indecomposable automorphic representation $V$ of ${\rm PSL}_2({\mathbb R})$ by lifting the weight 2 Eisenstein series $E_2$. Viewed as the Lie algebra $psl_2=psl_2({\mathbb R})$ representation, it can be identified with the Laurent polynomial ring
\begin{eqnarray}\nonumber
V\approx {\mathcal O}({\mathbb C}^\times)={\mathbb C}[z^{\pm 1}].
\end{eqnarray}
with the $psl_2$ action given by 
$$e={d\over{dz}},~~h=2z{d\over{dz}},~~f=z^2{d\over{dz}}.$$
The space $V$ has a natural multiplication map 
\begin{eqnarray}\label{8.19}\nonumber
V\to V\otimes V
\end{eqnarray}
which commutes with the
action of $psl_2$.  One can also identify the restricted dual space
with Laurent 1-forms
\begin{eqnarray}\nonumber
V'\approx \Omega^1({\mathbb C}^\times)={\mathbb C}[z^{\pm 1}]dz,
\end{eqnarray}
with the pairing of $V$ and $V'$ given by the residue of the product at 0.
This gives the dual map
\begin{eqnarray}\label{8.4}
V'\to V'\otimes V',
\end{eqnarray}
which can be used for the inductive limit construction as in
\cite{jones1,jones2,OS}.

\medskip

\noindent Namely for any polygon $P_\gamma$ with frontier $\gamma$, we assign a tensor
product
\begin{eqnarray}\nonumber
H_\gamma=\bigotimes_{j=1}^n V',
\end{eqnarray}
where $n$ is the number of edges in $\gamma$.  There is a partial ordering on the set of polygons given by inclusion $P_\gamma\subseteq P_{\gamma'}$, which induces the partial ordering
$\gamma\leq\gamma'$ on frontiers.  Clearly the Lie subalgebra 
$$ppsl_2(\gamma)\subseteq ppsl_2$$
of piecewise $psl_2$ maps with breakpoints at the vertices of $\gamma$ acts on
$H_\gamma$ in the natural way, where $psl_2$ corresponding to the $j^{\rm th}$ edge acts on the $j^{\rm th}$ factor of $H_\gamma$.  Then the map (\ref{8.4}) induces a family of intertwining operators
\begin{eqnarray}\label{8.6}\nonumber
T_\gamma^{\gamma'}:H_\gamma\to H_{\gamma'},
\end{eqnarray}
so that $T_\gamma^\gamma=1$ and the consistency condition $T^\gamma_{\gamma''}=
T^{\gamma'}_{\gamma''} T^\gamma_{\gamma'}$ holds whenever $\gamma\leq\gamma'\leq\gamma''$. The direct limit
\begin{eqnarray}
H=\varinjlim H_\gamma
\end{eqnarray}
yields a representation of $ppsl_2$. It can be made unitary by means of a 1-cocycle determined by the one-dimensional subrepresentation $V_0\subseteq V$ as in \cite{VGG} or alternatively
using a modification of the construction above based upon the approximation of the representation $V'$ of $psl_2$ by the representations of the complementary series $V(\lambda)$, for
$0<\lambda<2$ with $\lambda\neq1$, so that $V(\lambda)\approx V(2-\lambda)$ with
$\lim\limits_{\lambda\to 2} V(\lambda)=V$ and $\lim\limits_{\lambda\to 0} V(\lambda)=V'$;
see \cite{VGG} for details.

\begin{remark} The partial ordering above on ideal polygons occurring as fundamental domains for punctured surfaces leads one to consider
the so-called punctured solenoid $\mathcal H$, which is an initial object for the category of punctured surfaces with
morphisms given by finite covers branched only over the punctures.
The unpunctured version was introduced and studied by Dennis Sullivan in \cite{Sulsol}. 
 The decorated Teichm\"uller space of
$\mathcal H$ was studied in \cite{PS} and parametrized by certain coordinates on the Farey tesselation $\tau_*$, namely,
lambda length {\sl functions}, one such
continuous and ${\rm PSL}_2({\mathbb Z})$-equivariant function from the profinite completion of ${\rm PSL}_2({\mathbb Z})$
to ${\mathbb R}_{>0}$ for each edge of $\tau_*$. Flip-like generators for the mapping class group of $\mathcal H$ are derived in \cite{PS}, and a complete set of relations in these generators is finally given in \cite{BPS}.
\end{remark}

\noindent
To implement our program, we want to realize $H$, or heuristically $(V')^{\otimes\infty}$,
as an automorphic representation on ${\rm PPSL}_2({\mathbb Z})\backslash{\rm PPSL}_2({\mathbb R})$ of $ppsl_2$ .  To this end, we notice first that the space of harmonic functions
on ${\rm PSL}_2({\mathbb R})$, cf.~(\ref{8.18}), can be identified with the space of functions on the boundary of ${\rm PSL}_2({\mathbb R})$, which is dual to the space of functions of the holographic coordinates $\{ (s,\delta)\}$ introduced in Section 2.  Furthermore for both spaces, we can also identify the corresponding ${\rm PSL}_2({\mathbb R})$ actions, and therefore, we can
realize $V'$ by imposing modular invariance on the space of functions depending on a single pair $(s,\delta)$ of holographic coordinates.
This realization clearly extends to finite tensor products $(V')^{\otimes n}$ associated to $n$-gons $P_\gamma$.  We conjecture that the inductive limit $H$ also admits a realization by the automorphic functions in holographic coordinates, with ${\rm PPSL}_2({\mathbb Z})$ playing the role of the modular group ${\rm PSL}_2({\mathbb Z})$.  

\medskip

\noindent
On the other hand, the analogy between ${\rm PPSL}_2({\mathbb Z})$ and the mapping class group
suggests an extension of the former by an infinite symmetric group, or more generally by an infinite braid group, which is still a discrete subgroup of ${\rm Homeo}_+({\mathbb S}^1)$,
cf.~\cite{FK} and the references therein.  The reduction of $H$ by the additional symmetry yields the space $\Lambda$, or heuristically $S^\infty V'$, which is naturally identified with the free bsosonic field, one of the simplest ${\rm CFT}_2$.

\medskip

\noindent 
As discussed in Appendix B, the theory of automorphic forms allows the enlargement of various spaces by relaxing the automorphic property to smaller discrete subgroups; for example, the restriction to invariance by the commutant ${\rm PSL}_2({\mathbb Z})'\subseteq {\rm PSL}_2({\mathbb Z})$ enlarges $V$ sixfold additively. Similarly, we expect that replacing
${\rm PPSL}_2({\mathbb Z})$ by ${\rm P}({\rm PSL}_2({\mathbb Z}))'$ will
correspondingly yield a multiplicative increase of $\Lambda$ and thereby allow the construction of
a ${\rm CFT}_2$ of size comparable to the Monster ${\rm CFT}_2$.


\medskip

\noindent 
Finally we come to the question of how we can effectively capture the Monster using the new universal automorphic theory, and in particular how we can derive  the Monstrous Moonshine.
At this moment, we do not know the complete answer, but we can begin on the opposite side, which will help to clarify the problem and capture the beast.  Namely, we consider all 194 conjugacy classes of the Monster and the Thompson series of the representations of these classes in the the Monster representation of \cite{FLM}.  Then we know from \cite{B} that they are the canonical Hauptmoduln $J_\Gamma(z)$ for some genus-zero discrete subgroups $\Gamma$ of ${\rm PSL}_2({\mathbb R})$ of Moonshine type, i.e., groups commensurable with ${\rm PSL}_2({\mathbb Z})$ which contain $\Gamma_\infty=\{\pm\begin{psmallmatrix}1&n\\0&1\end{psmallmatrix}\}$.  Then one has
\begin{eqnarray}\nonumber
J_\Gamma (z)=q^{-1}+\sum_{n=1}^\infty c_\Gamma(n)q^n,
\end{eqnarray}
for some integral coefficients $c_\Gamma(n)$.
One can also construct the corresponding Eisenstein series using a logarithmic derivative
of $J_\Gamma(z)$, or as usual an analytic continuation to $s=0$ of the series
\begin{eqnarray}\nonumber
E_2^\Gamma(z,s)=\sum_{(c,d)\in\Gamma/\Gamma_\infty}
{{({\rm Im}~z)^s}\over{(cz+d)^2|cz+d|^{2s}}}.
\end{eqnarray}
As in the classical example $\Gamma={\rm PSL}_2({\mathbb Z})$, the resulting
Eisenstein series after multiplication by $dz$ is $\Gamma$-invariant though not holomorphic.
Then we can repeat the lift of a $\Gamma$-invariant 1-form recalled in Appendix B to the automorphic function
\begin{eqnarray}\nonumber
f=E_2^\Gamma\rightsquigarrow \phi_f(g)=(ci+d)^{-2}f(g\cdot i),
\end{eqnarray}
where $g=\begin{psmallmatrix}a&b\\c&d\\\end{psmallmatrix}\in{\rm PSL}_2({\mathbb R})$,
and we again obtain an indecomposable representation of ${\rm PSL}_2({\mathbb R})$ that we shall denote $V^\Gamma$.  In fact, we just get another model of the same indecomposable
representation of ${\rm PSL}_2({\mathbb R})$.  We can also expect a similar construction of the space of universal automorphic forms for ${\rm P}\Gamma$ denoted $\Lambda^\Gamma$.  It is natural to conjecture that we find a twisted conformal field theory associated to $\Lambda^\Gamma$ and the Monster conjugacy class corresponding to $\Gamma$.  The comparison of $\Lambda$ and $\Lambda^\Gamma$ should yield the Monster element in this conjugacy class, and they together will allow realization of the Monster group via automorphisms of $\Lambda$.  At this point, the 
trap will close upon the Monster, and the Moonshine will follow.

\medskip


\bigskip

\appendix

\section {Harmonic analysis}\label{appendixA}

\bigskip

\noindent To begin, we express the generators
of the Witt algebra in terms
of the normalized arithmetic wavelets as computed in Theorem 4.13 of \cite{Pbook}
and going back in essence to \cite{MP}.

\begin{theorem}\label{thm:Philbert}
For each $n\in{\mathbb Z}$, we have the wavelet expansion
$$\begin{aligned}
L_n=e^{in\theta}{\partial\over{\partial\theta}}&=(b_0^n+b_{+1}^ne^{i\theta}+b_{-1}^ne^{-i\theta}){\partial\over{\partial\theta}}\\
&~~+~{i\over 4}~\sum _{e\in\tau _*}\biggl\{ n(\xi ^n+\eta ^n)+{{\eta +\xi}\over{\eta
-\xi}}~(\xi ^n-\eta ^n)\biggr\}~{\overline{\vartheta}}_e(\theta),
\end{aligned}$$
where $e\in\tau _*$ has
ideal points $\xi ,\eta \in {\mathbb S}^1$ and
$$b_0^n=\begin{cases} +1, n\equiv 0(4);\\ \hskip 1.5ex 0, n\equiv 1(4);\\ +1, n\equiv 2(4);\\ \hskip 1.5ex 0, n\equiv 3(4);\\ \end{cases} b_{+1}^n=\begin{cases} \hskip 1.5ex 0, n\equiv 0(4);\\ +1, n\equiv 1(4);\\ -i, n\equiv 2(4);\\ \hskip 1.5ex 0, n\equiv 3(4);\\ \end{cases} b_{-1}^n=\begin{cases} \hskip 1.5ex 0, n\equiv 0(4);\\ \hskip 1.5ex 0, n\equiv 1(4);\\ +i, n\equiv 2(4);\\ +1, n\equiv 3(4).\\ \end{cases} $$
\end{theorem}

\medskip

\noindent This is not the simplest expression $e^{in\theta}\doteq\sum_{|e|\in\tau_*} g_n^e \bar\vartheta_{|e|}$, rather these particular coefficients $g_n^e$ are specially chosen in Theorem~\ref{thm:Philbert} to guarantee their suitable decay in $n$, cf.\ Theorem 6.4 of \cite{Philbert}.

\medskip

\noindent The next result, Theorem 4.11 of \cite{Pbook}, which goes back to \cite{Puniv}, is of basic utility and provides the Fourier expansion
of the normalized arithmetic wavelets.

\begin{theorem}\label{thm:fourier}
If $A=\begin{psmallmatrix}a&b\\c&d\end{psmallmatrix}\in {\rm PSL}_2
$, then the Fourier expansion $\bar\vartheta_A\sim\sum _{n\in{\mathbb Z}} c_n~e^{in\theta}{\partial\over{\partial\theta}}$  for $n^2>1$ is given by
$$\aligned
\pi i(n^3-n)~~c_n&~~=~~~-[(c-a)^2+(b-d)^2]\biggl [{{(b-d)-i(a-c)}\over
{(b-d)+i(a-c)}}\biggr ] ^n\\
&~~~~~~~~+2(c^2+d^2)\biggl [ {{d-ic}\over{d+ic}}\biggr ]^n
~~+~~2(a^2+b^2)\biggl [ {{b-ia}\over{b+ia}}\biggr ]^n\\
&~~~~~~~~-[(c+a)^2+(b+d)^2]\biggl [{{(b+d)-i(a+c)}\over
{(b+d)+i(a+c)}}\biggr ] ^n,
\endaligned$$
where the Fourier modes $c_0,c_{\pm 1}$ are chosen to guarantee that the expansion is normalized.
 In particular for the mother wavelet, we have
$\bar\vartheta=\vartheta \sim  {8\over{\pi i}}~\sum _{n\equiv 2(4)}~{1\over{n^3-n}}~e^{in\theta}{\partial\over{\partial\theta}}$.
\end{theorem}

\begin{remark} \label{rmkonKK}
A sketch of the proof is that since $\bar\vartheta_A$ is once continuously
differentiable, we can twice integrate by parts the standard expression for the Fourier
coefficients $c_n$ with $n^2>1$ and derive these expressions without difficulty.   For $n^2\leq 1$ the corresponding
equations give certain constraining equations on the zeroeth Fourier modes,
which are crucial for the computation of the Kirillov-Kostant
form, cf. Section \ref{sec:2-forms}
\end{remark}

\noindent It is not difficult to compute directly the Fourier expansion of the (unnormalized) hyperfans as follows.

\begin{proposition}\label{hyperfanfourier}
Given $A=\begin{psmallmatrix}a&b\\c&d\\\end{psmallmatrix}\in{\rm PSL}_2$,
let $$\zeta_-={{b-ia}\over{b+ia}},~~\zeta_+={{d-ic}\over{d+ic}}\in {\mathbb S}^1.$$
Then the hyperfan $\psi_A$ 
has Fourier expansion $\psi_A\sim$ $\sum c_n~e^{in\theta}$ with
$$\begin{aligned}2\pi i ~c_n&=\zeta_-^n\biggl[
{{(d-ic)^2}\over{2(n+1)}}\zeta_-+{{(d+ic)^2}\over{2(n-1)}}\zeta_-^{-1}-{{c^2+d^2}\over n}\biggr]\\
&-\zeta_+^n\biggl[
{{(d-ic)^2}\over{2(n+1)}}\zeta_++{{(d+ic)^2}\over{2(n-1)}}\zeta_+^{-1}-{{c^2+d^2}\over n}\biggr],\\
\end{aligned}$$
for $n^2>2$, and
$$\begin{aligned}
2\pi i~c_0&={{(d-ic)^2}\over 2}[\zeta_+^{-1}-\zeta_-^{-1}]-{{(d+ic)^2}\over 2}[\zeta_+-\zeta_-]\\
&\hskip 2.7ex+i(c^2+d^2)[\theta_+-\theta_-],\\
2\pi i~c_{\pm1}&=\mp{{(d\pm ic)^2}\over 2}[\zeta_+^{\pm 2}-\zeta_-^{\pm2}]\pm(c^2+d^2)[\zeta_+^{\pm1}-\zeta_-^{\pm1}]\\
&\hskip 2.7ex -{{i(d\mp ic)^2}\over 2}[\theta_+-\theta_-],\\
\end{aligned}$$
where $\theta_\mp={\rm arctan}\, \zeta_\pm$.
\end{proposition}

\begin{proof}
$\psi_A$ is supported on an interval with respective initial and terminal points given by 
$\zeta_-^{-1}$ and $\zeta_+^{-1}$, and takes there the values
$$A^{-1}\begin{pmatrix}0&1\\0&0\end{pmatrix}A=\begin{pmatrix}
\hskip .5excd&\hskip 2exd^2\\-c^2&-cd\end{pmatrix}$$
or in other words 
${{(d-ic)^2}\over 2}\,e^{i\theta}+{{(d+ic)^2}\over 2}\,e^{-i\theta}-(c^2+d^2)$.  The standard expression for Fourier coefficients yields the results.
\end{proof}

\noindent Notice that the normalized hyperfans $\bar\psi_A$ have the same higher Fourier coefficients as $\psi_A$ since $\psi_A\doteq\bar\psi_A$.  The small Fourier modes
for $\bar\psi_A$, and likewise those of $\bar\theta_A$, can be computed with some difficulty by calculating the values of these vector fields at ${0\over 1},{1\over 1},{1\over 0}\in{\mathbb S}^1$.  As we do not require them, we omit the details.

\bigskip

\section{The $E_2$ automorphic representation}

\bigskip

\noindent The automorphic form that we have constructed on the universal Teichm\"uller space has a simple analogue on the hyperbolic plane, namely,
\begin{eqnarray}\nonumber
\xi(z)=E_2(z)dz,
\end{eqnarray}
where $E_2(z)$ is an almost holomorphic modular form of weight 2.  Since $dz$ transforms with respect to ${\rm PSL}_2({\mathbb R})$ with the factor balancing that of $E_2(z)$, we obtain a 
${\rm PSL}_2({\mathbb Z})$-invariant 1-form that can be viewed as the classical counterpart of the automorphic form constructed in Section 7.  

\medskip

\noindent To state our results unambiguously, let us recall the following explicit formulae, c.f.\,\cite{DS}.  Set
\begin{eqnarray}\nonumber
{\mathbb E}_2(z)=1-24\sum_{n=1}^\infty\sigma(n)q^n,
\end{eqnarray}
where $\sigma(n)$ is the sum of all the positive divisors of $n$ and we set $q={\rm exp}(2\pi i z).$
Then ${\mathbb E}_2(z)$ is holomorphic but not quite modular invariant since
\begin{eqnarray}\nonumber
{\mathbb E}_2(-{1\over z})=z^2{\mathbb E}_2(z)+{{12}\over{2\pi i}}z.
\end{eqnarray}
To get an honest weight 2 form, we need a correction
\begin{eqnarray}\nonumber
E(z)={\mathbb E}_2(z)-{3\over{\pi~{\rm Im}(z)}},
\end{eqnarray}
but we thereby lose the holomorphicity of ${\mathbb E}_2(z)$, since
\begin{eqnarray}\nonumber
{d\over{d\bar{z}}}E(z)={{3i}\over{2\pi}}{1\over{({\rm Im}(z))^2}}.
\end{eqnarray}
In the classical theory of automorphic forms, one can lift the holomorphic modular forms from the hyperbolic plane to ${\rm PSL}_2({\mathbb R})$ so that they become the lowest weight vectors of the holomorphic series representations.  Specifically for a weight $2k$ holomorphic form
$f:{\mathcal U}\to{\mathbb C}$, where $k\geq 2$, we define
\begin{eqnarray}\nonumber
\phi_f(g)=(ci+d)^{-2k} f(g\cdot i),
\end{eqnarray}
as in \cite{GGP}, where $g=\begin{psmallmatrix}1&b\\c&d\\\end{psmallmatrix}\in{\rm PSL}_2({\mathbb R})$.
It follows that $\phi_f$ is invariant under the left action 
\begin{eqnarray}\label{8.7}
\phi_f(\gamma g)=\phi_f(g),~{\rm for}~ \gamma\in{\rm PSL}_2({\mathbb Z})
\end{eqnarray}
of ${\rm PSL}_2({\mathbb Z})$. One can then define the action of ${\rm PSL}_2({\mathbb R})$ on a certain space of automorphic functions generated by $\phi_f$ by the right action
\begin{eqnarray}\nonumber
(\pi(g)\phi)(h)=\phi(hg),~{\rm for}~h,g\in{\rm PSL}_2({\mathbb R}),
\end{eqnarray}
so that $\phi_f$ becomes the lowest weight vector for the Lie algebra
$psl_2({\mathbb R})$ in the appropriate basis $E,F,H$.  

\medskip

\noindent To state the explicit
formulae, we shall rely on explicit calculations from \cite{FGKP}.  Let
$h,e,f$ be the standard basis of $psl_2({\mathbb R})$ and consider the Cayley
transform of this basis
\begin{eqnarray}\nonumber
H=\begin{pmatrix}0&-1\\i&~~0\end{pmatrix},
E={1\over 2}\begin{pmatrix}1&~~i\\i&-1\end{pmatrix},
F={1\over 2}\begin{pmatrix}~~1&-i\\-i&~~1\end{pmatrix}.
\end{eqnarray}
Using the standard parametrization 
\begin{eqnarray}\label{8.11}
\begin{pmatrix}a&b\\c&d\\\end{pmatrix}=
\begin{pmatrix}1&x\\0&1\\\end{pmatrix}
\begin{pmatrix}y^{1\over 2}&0\\0&y^{-{1\over 2}}\\\end{pmatrix}
\begin{pmatrix}~~{\rm cos}\,\theta&{\rm sin}\,\theta\\-{\rm sin}\,\theta&{\rm cos}\,\theta\\\end{pmatrix}
\end{eqnarray}
of ${\rm PSL}_2({\mathbb R})$,
one can find an explicit action of the Lie algebra basis as follows
\begin{eqnarray}\nonumber
H&=&-\partial_\theta,\\
E&=&2i e^{2i\theta}(y\partial_{\bar{z}}-{1\over 4}\partial_\theta),\nonumber\\
F&=&-2i e^{-2i\theta}(y\partial_{\bar{z}}-{1\over 4}\partial_\theta).\nonumber
\end{eqnarray}
Thus, one has
\begin{eqnarray}\nonumber
F\phi_f&=&0,\\\nonumber
H\phi_f&=&2k\phi_f.
\end{eqnarray}

\medskip

\noindent
In the case that $k=1$, the first of these two equations no longer holds, and instead one has
\begin{eqnarray}\nonumber
F\phi_f={3\over \pi}.
\end{eqnarray}
It follows that $\phi_f$ is no longer a highest weight vector, but it also generates
a one-dimensional sub-representation.  Only after factorization with respect to this
sub-representation does one obtain the weight 2 irreducible representation of the holomorphic
series, denoted $V_+$.

\medskip

\noindent
Similarly using the complex conjugation, one gets the weight 2 representation of the antiholomorphic series denoted $V_-$.  The one-dimensional sub-representation, denoted $V_0$, is common to both the weight 2 indecomposable representations, and we get a larger indecomposable representation, denoted $V$, which has a simple composition series
\begin{eqnarray}\nonumber
0\to V_0\to V\to V_+\otimes V_-\to 0.
\end{eqnarray}

\medskip

\noindent
In order to characterize the resulting three component representation, one can use the Casimir
operator, which in the explicit coordinates (\ref{8.11}) has the form
\begin{eqnarray}\nonumber
\Delta=y^2({{\partial^2}\over{\partial x^2}}+{{\partial^2}\over{\partial y^2}})-y{{\partial^2}\over{\partial x\partial\theta}}.
\end{eqnarray}
It annihilates the weight 2 automorphic representation generated by $\Phi_f$ and its conjugate, where $f=E_2$, to wit
\begin{eqnarray}\label{8.18}
\Delta \phi_{E_2}=0.
\end{eqnarray}
One can check that this harmonic property (\ref{8.18}) together with the automorphicity condition
(\ref{8.7}) and the usual growth condition of automorphic functions characterizes the three component representation $V$ that we have obtained by lifting the $E_2$ modular 1-form.
The problem of fundamental importance in our new theory is to construct a representation of the Lie algebra $ppsl_2({\mathbb R})$ using a lift of our automorphic 1-form in Section 7.

\medskip

\noindent
The automorphic representation $V$ that we have constructed is well known in the representation theory of ${\rm PSL}_2({\mathbb R})$ as a limit of complementary series \cite{VGG}.  The representations $V_+$ and $V_-$ are the spin 1 irreducible representations of the holomorphic and antiholomorphic discrete series.  Note that the lift of Eisenstein series $E_{2k}$, for $k=2,3,\ldots$, yields spin $k$ irreducible representations of the same series.  The representation $V$ plays a pivotal role and somehow seems to have been missed in the theory of automorphic functions for the pair ${\rm PSL}_2({\mathbb Z})\subseteq{\rm PSL}_2({\mathbb R})$
Thus it is expected that the counterpart of the representation $V$ in our new theory will play an equally fundamental role, and it is one of our main open problems to construct it explicitly.

\medskip 

\noindent
One can ask what is an advantage of a realization of this representation as an automorphic representation.  One such benefit is that we can naturally increase this representation by relaxing
the automorphicity condition (\ref{8.7}) from the modular group ${\rm PSL}_2({\mathbb Z})$ 
to various subgroups including its commutant ${\rm PSL}_2({\mathbb Z})'$ considered in Section 1.  This immediately yields six copies of the weight 2 representation.

\medskip

\noindent
Finally, we note that we could also consider the weight 1 representations which correspond to the limit of the holomorphic discrete series of ${\rm SL}_2({\mathbb R})$.  The index
of the commutant  is  $[{\rm SL}_2({\mathbb R}):{\rm SL}_2({\mathbb R})']=12$ in this case, and we get correspondingly 12 copies upon relaxing the automorphicity condition.


\begin{thebibliography}{ABCD}

\bibitem{B} L.\ Bers, ``Universal Teichm\"uller space", in {\it Analytic Methods in Mathematical Physics}, pp. 65-83, Gordon and Breach, New York, 1968.

\bibitem{BPS} S.\ Bonnot, R.\ Penner, D.\ Saric, ``A presentation for the baseleaf preserving mapping class group of the punctured solenoid", {\it Algebraic and Geometric Topology} {\bf 7}(3) (2007), 1171-1199.

\bibitem{Milgram}
C.\ P.\ Boyer, J.\ C.\ Hurtubise, B.\ M.\ Mann, R.\ J.\ Milgram, ``The topology 
of instanton moduli spaces, I: the Atiyah-Jones Conjecture", {\it Annals of Mathematics} {\bf 137} (1993), 561-609.

\bibitem{Brown} F.\ Brown, ``Multiple zeta values and periods of moduli spaces $M_0^n$", {\it Annalens Scientifiques de l'\'Ecole Normale Sup\'erieure} {\bf 42} (2009), 371-489.

\bibitem{CFP} J.\ Cannon, W.\ Floyd, W.\ Parry, ``Introductory remarks on Richard Thompson's groups", L'Enseignement Mathematique {\bf 42} (1996), 215-256.

\bibitem{171} J.\ Conway, J.\ McKay, A.\ Sabbar, ``On the discrete groups of Moonshine", {\it Proceedings of the American Mathematical Society} {\bf 132} (8) (2004), 2233-2240.

\bibitem{CN} J.\ Conway and S.\ Norton, ``Monstrous moonshine", {\it Bulletin London Mathematical Society} {\bf 11} (1979), 308-339.

\bibitem{C} C.\ Cummins, ``Congruence subgroups of groups commensurable with ${\rm PSL}_2({\mathbb Z})$ of genus 0 and 1", {\it Experimental Mathematics} {\bf 13} (2004), 361-382.

\bibitem{DS} F.\ Diamond and J.\ Shurman, {\it A first course in modular forms},
Springer-Verlag, 2005.

\bibitem{DF} J.\ Duncan and I.\ Frenkel, ``Rademacher sums, moonshine and gravity'',
{\it Communications in Number Theory and Physics} {\bf 5}(4) (2011), 1-128.s"

\bibitem{FGKP} P.\ Fleig, H.\ Gustaffson, A.\ Kleinschmidt, D.\ Persson,
``Eisenstein series and automorphic representations", {\it Cambridge Studies in Advanced Mathematics} Cambridge University Press, 2018.

\bibitem{FLM} I.\ Frenkel, J.\ Lepowsky, A.\ Meurman, ``Vertex operator algebras and the Monster", {\it Pure and Applied Mathemartics} {\bf 134}, Academic Press, Boston, MA 1988.

\bibitem{Fuchs}
J.\ Fuchs, Affine Lie Algebras and Quantum Groups, Cambridge University Press (1992).

\bibitem{FK} L.\ Funar and C.\ Kapoudjian, ``The braided Ptolemy-Thompson group is finitely presented'', {\it Geometry and Topology} {\bf 12} (2008), 475-530.


\bibitem{GGP} I.\ Gelfand, M.\ Graev and I.\ Piatetski-Shapiro, {\it Representation theory
and automorphic functions}, Saunders Mathematics Books, Saunders, 1968.




\bibitem{jones1} V.\ F.\ R.\ Jones. ``Some unitary representations of Thompson?s groups F and T", {\it Journal of Combinatorial Algebra} {\bf 1}(1) (2014), 1-44.
arXiv: 1412.7740.

\bibitem{jones2} V.\ F.\ R.\ Jones. ``A no-go theorem for the continuum limit of a periodic quantum spin chain", {\it Communications in Mathematical Physics} {\bf 357} (2018), 295-317.

\bibitem{KK}
A.\ A.\ Kirillov, Lectures on the Orbit Method, Graduate Studies in Mathematics, Vol. 64, American Mathematical Society.

\bibitem{L} G.\ Laget, ``Groupes de Thompson projetifs de genre 0", These de doctorat
en Mathematiques, Universite Joseph-Fourier, Grenoble, France 2004.

\bibitem{MP} F.\ Malikov and R.\ C.\ Penner, ``The Lie algebra of homeomorphisms
of the circle", {\it Advances in Mathematics} {\bf 140} (1999), 282-322.

\bibitem{Mink}
H.\ Minkowski,  "Zur Geometrie der Zahlen", Verhandlungen des III. internationalen Mathematiker-Kongresses in Heidelberg, Berlin (1907), pp. 164-173, JFM 36.0281.01, archived from the original on 4 January 2015.

\bibitem{OS} T.\ J.\ Osborne and D.\ E.\ Steigemann, ``Dynamics for holographic codes",
arXiv:1706.08823.

\bibitem{Pdec} R.\ C.\ Penner, ``The decorated Teichm\"uller space of punctured surfaces'', 
{\it Communications in Mathematical Physics} {\bf 113} (1987), 299-339.

\bibitem {Puniv} ---, ``Universal constructions in Teichm\"uller theory'', 
{\it Advances in Mathematics} {\bf 98} (1993), 143-215.

\bibitem{PWP} ---, ``Weil-Petersson volumes'', {\it Journal of Differential Geometry}
{\bf 35} (1992), 559-608.

\bibitem{Pbook} ---, {\it Decorated Teichm\"uller theory},
QGM Masters Class Series, volume 1, European
Math Society  (2012).

\bibitem{Philbert} ---, ``On Hilbert, Fourier and wavelet transforms'',
{\it Communications on Pure and Applied Mathematics} {\bf 15} (2002), 772-814.

\bibitem{PS} --- and D.\ S\'aric, ``Decorated Teichm\"uller theory of the punctured solenoid'',
{\it Geometriae Dedicata} {\bf 132} (2008), 179-212.


\bibitem{SK}
C.\ Scainci and K.\ Krasnov, ``The Universal Phase Space of AdS$_3$ Gravity", {\it Communications in Mathematical Physics} {\bf 322} (2013), 167-205.

\bibitem{Segal}
A.\ Pressley and G.\ Segal, Loop Groups, Oxford Mathematical Monographs (1986), Oxford University Press.


\bibitem{Sulsol}
D.\ Sullivan, ``Linking the universalities of Milnor-Thurston, Feigenbaum and AhlforsBers", Milnor Festschrift, Topological methods in modern mathematics (L.\ Goldberg and
A.\ Phillips, eds.), Publish or Perish (1993), 543-563.


\bibitem{Verjovsky-Nag} A.\ Verjovsky and S.\ Nag, ``Diff($S^1$) and the Teichm\"uller spaces", {\it Communications in Mathematical Physics} {\bf 130} (1990), 123-138.

\bibitem{V} H.\ Verlinde, ``Conformal field theory, two-dimensional quantum gravity and quantization of Teichm\"uller space", {\it Nuclear Physics B} {\bf 337} (1990), 652-680.

\bibitem{VGG} A.\ Vershik, I.\ Gelfand, M.\ Graev, ``Representations of the group $SL(2,{R})$, where $R$ is a ring of functions", {\it Uspekhi Mat. Nauk} {\bf 28} (1973), 83-128.

\bibitem{W} E.\ Witten, ``Three-dimensional quantum gravity revisited", arXiv:0706.3359, June 2007.

\end{thebibliography}
\end{document}